\documentclass[a4paper,11pt]{amsart}
\usepackage{amsfonts}
\usepackage{amsmath}
\usepackage{amssymb}
\usepackage{dutchcal}
\usepackage{float}
\usepackage{mathrsfs}
\usepackage{geometry}
\usepackage{graphicx}
\usepackage[pdftex,colorlinks=true,linkcolor=blue,citecolor=lessblue,breaklinks=true,bookmarksopen=true]{hyperref}
\usepackage[utf8]{inputenc}
\usepackage{leftidx} 
\usepackage{mathrsfs}
\usepackage{mathtools}
\usepackage{pdflscape}
\usepackage{tikz}
\usepackage{xcolor}
\definecolor{winered}{rgb}{0.7,0,0}
\definecolor{lessblue}{rgb}{0,0,0.7}
\usepackage{xspace} 
\DeclareMathSymbol{\intop}{\mathop}{symbols}{115}
\def\int{\intop}

\DeclareMathSymbol{\sumop}{\mathop}{largesymbols}{80}
\def\ssum{\sumop}
\let \sum \relax
\DeclareMathOperator*{\sum}{\textstyle{\ssum}}

\DeclareMathSymbol{\prodop}{\mathop}{largesymbols}{81}
\def\sprod{\prodop}
\let \prod \relax
\DeclareMathOperator*{\prod}{\textstyle{\sprod}}

\DeclareMathSymbol{\bigoplusop}{\mathop}{largesymbols}{76}

\newcommand{\mb}[1]{\ensuremath{\mathbb{#1}}}
\newcommand{\N}{{\mb{N}}}

\newcommand{\R}{{\mb{R}}}
\newcommand{\C}{{\mb{C}}}

\newcommand{\eps}{\varepsilon}

\newcommand{\B}{{\mathscr B}}
\newcommand{\D}{{\mathscr D}}

\renewcommand{\H}{{\mathcal H}}

\newcommand{\caO}{\mathcal O}

\newcommand{\M}{{\mathcal M}}

\newcommand{\U}{\ensuremath{\mathscr U}}
\newcommand{\caU}{\ensuremath{\mathcal U}}

\let \Re \relax
\DeclareMathOperator{\Re}{Re}
\let \Im \relax
\DeclareMathOperator{\Im}{Im}

\let \div \relax
\DeclareMathOperator{\div}{div}

\newcommand{\ovl}[1]{\overline{#1}}

\newcommand{\Con}{\ensuremath{\mathscr C}}

\newcommand{\Cinf}{\ensuremath{\Con^\infty}}

\newcommand{\Cbarc}{\ensuremath{\overline \Con^\infty}_c}

\newcommand{\Cinfc}{\ensuremath{\Con^\infty}_{c}}
\renewcommand{\S}{\ensuremath{\mathscr S}}
\renewcommand{\D}{\ensuremath{\mathscr D}}

\newcommand{\inp}[2]{( #1, #2 )} 
\newcommand{\biginp}[2]{\big( #1, #2 \big)}

\newcommand{\Op}{\ensuremath{\mathrm{Op}}}

\newcommand{\Opt}{\ensuremath{\mathrm{Op}_{\mbox{\tiny ${\mathsf T}$}}}}
\newcommand{\opt}{\ensuremath{\mathrm{Op}_{\mbox{\tiny ${\mathsf T}$}}}}
\DeclareMathOperator{\id}{Id}
\DeclareMathOperator{\supp}{supp}

\renewcommand{\d}{\ensuremath{\partial}}

\renewcommand{\qedsymbol}{$\blacksquare$}

\newcommand{\bld}[1]{\mbox{\boldmath $#1$}}

\newcommand{\rhs}{r.h.s.\@\xspace}
\newcommand{\lhs}{l.h.s.\@\xspace}

\newcommand{\eg}{e.g.\@\xspace}

\newcommand{\se}{sub-ellipticity\xspace}
\newcommand{\LS}{Lopatinski\u\i-\v{S}apiro\xspace}
\newcommand{\nhd}{neighborhood\xspace}

\newcommand{\st}{such that\xspace}
\newcommand{\suff}{sufficiently\xspace}

\newcommand{\br}{\ensuremath{_{|x_d=0^+}}}
\newcommand{\bt}{\ensuremath{_{|{x_d}=0}}}

\newcommand{\ic}{\ensuremath{_{|t=0}}}
\newcommand{\Ssc}{S_\tau}

\newcommand{\Psisc}{\Psi_\tau}
\newcommand{\Psiscl}{\Psi_{\tau,\mathrm{ph}}}
\newcommand{\Dsc}{\D_\tau}

\newcommand{\Ssct}{S_{\mbox{\tiny ${\mathsf T}$},\tau}}

\newcommand{\Psisct}{\Psi_{\mbox{\tiny ${\mathsf T}$},\tau}}

\newcommand{\Dsct}{\D_{\mbox{\tiny ${\mathsf T}$},\tau}}

\newcommand{\Norm}[2]{{\| #1 \|}_{#2}}

\newcommand{\norm}[2]{{|#1|}_{#2}}

\newcommand{\Normsc}[2]{{\| #1 \|}_{\tau,#2}}
\newcommand{\normsc}[2]{{| #1 |}_{\tau,#2}}

\newcommand{\lsc}{\lambda_\tau}
\newcommand{\lscbis}[1]{\lambda_{\tau, #1}}
\newcommand{\lsct}{\lambda_{\mbox{\tiny ${\mathsf T}$},\tau}}
\newcommand{\lsctbis}[1]{\lambda_{\mbox{\tiny ${\mathsf T}$},\tau, #1}}

\newcommand{\Lsc}{\Lambda_\tau}
\newcommand{\Lsct}{\Lambda_{\mbox{\tiny ${\mathsf T}$},\tau}}

\newcommand{\tchi}{\tilde{\chi}}

\newcommand{\tpi}{\tilde{\pi}}
\newcommand{\tbeta}{\tilde{\beta}}

\newcommand{\ttau}{\tilde{\tau}}
\newcommand{\csp}{\gamma}

\newcounter{sectionapp}

\DeclareMathOperator{\trace}{tr}

\newcommand{\Rdp}{\R^d_+}

\newcommand{\energy}{\mathcal E}

\newcommand{\nablag}{\nabla_{\!\! g}}

\newcommand{\tW}{\tilde{W}}

\definecolor{lightpink}{rgb}{1,0.75,0.75}

\newcommand{\dnu}{\d_{\nu}}

\newcommand{\et}{\ensuremath{\text{and}}}

\newcommand{\avec}{\ensuremath{\text{with}}}
\newcommand{\si}{\ensuremath{\text{if}}}

\newcommand{\dans}{\ensuremath{\text{in}}}
\newcommand{\sur}{\ensuremath{\text{on}}}
\newcommand{\where}{\ensuremath{\text{where}}}

\newcommand{\Ellop}{P}
\newcommand{\ellop}{p}
\newcommand{\tellop}{\tilde{p}}
\newcommand{\Psig}{\Ellop_\sigma}
\newcommand{\psig}{\ellop_\sigma}

\newcommand{\pphi}{\ellop^\varphi}
\newcommand{\Psigphi}{\Psig^\varphi}
\newcommand{\psigphi}{\psig^\varphi}
\newcommand{\ppsi}{\ellop^\psi}

\newcommand{\Ddphi}{D_d^\varphi}

\newcommand{\SgGen}{A}

\newcommand{\imply}{\Rightarrow}

\newcommand{\Bsig}{B_\sigma}
\newcommand{\Bphi}{B^\varphi}
\newcommand{\bphi}{b^\varphi}

\newcommand{\gsf}{\mathsf g}
\newtheorem{theorem}{Theorem}[section]
\newtheorem{proposition}[theorem]{Proposition}
\newtheorem{lemma}[theorem]{Lemma}
\newtheorem{corollary}[theorem]{Corollary}

\theoremstyle{definition}
\newtheorem{definition}[theorem]{Definition}

\newtheorem{remark}[theorem]{Remark}

\numberwithin{equation}{section}
\makeatletter
\@namedef{subjclassname@2010}{\textup{2020} Mathematics Subject Classification}
\makeatother
\begin{document}
	
\pagenumbering{arabic}	
\title[Dynamical Boundary Condition]{Damped Wave Equation with a Dynamic Boundary Condition: Sharp Resolvent Estimates}

\author[Buffe, Gonzalez~Martinez, Le~Rousseau, and Robbiano]{R. Buffe, V. H. Gonzalez Martinez, J. Le Rousseau, and L. Robbiano}


\address{R\'emi Buffe. Universit\'e de Lorraine, CNRS, Inria, IECL, F-54000 Nancy, France}
\email{remi.buffe@univ-lorraine.fr}
\address{Victor Hugo Gonzalez Martinez. Departamento de Matem\'atica,  Universidade Federal de Pernambuco (UFPE), 50740-545, Recife (PE), Brazil.}
\email{victor.martinez@ufpe.br}
\address{Laboratoire Analyse, Géométrie et Applications, CNRS UMR 7539, Université Sorbonne Paris-Nord, F-93430 Villetaneuse, France.}
\email{gonzalezmartinez@math.univ-paris13.fr}

\thanks{*Corresponding author: jlr@math.univ-paris13.fr}
\address{Jérôme Le Rousseau. Laboratoire Analyse, Géométrie et Applications, CNRS UMR 7539, Université Sorbonne Paris-Nord, F-93430 Villetaneuse, France.}
\email{jlr@math.univ-paris13.fr}
\address{Luc Robbiano. Laboratoire de Mathématiques de Versailles, UVSQ, CNRS, Université Paris-Saclay, F-78035 Versailles, France.}
\email{luc.robbiano@uvsq.fr} 

\date{\today}
	
\begin{abstract}
This work investigates the stabilization of a wave equation
with a localized damping and dynamic boundary conditions. This model
traces back to the classical treatment of torsion and impact
problems. Consider a second-order hyperbolic equation posed with a
dynamic boundary condition on a smooth compact Riemannian manifold
$(\M,g)$ with boundary, in which the damping acts either in the
interior or on $\d \M$. Under the assumption that the damping is
positive on an open subset of its domain of action, a resolvent
estimate is established. A logarithmic decay rate is then deduced for
the energy of strong solutions. We provide examples that establish the
sharpness of the resolvent estimate.
\end{abstract}
\maketitle

\def\proofname{\bfseries Proof}

\tableofcontents

\section{Introduction}

The wave equation $\d_t^2 u - \Delta u = 0$ has been the subject of extensive and
sustained study, driven both by its rich mathematical structure and its
broad applicability in science and engineering. Among the many
variants of this equation that arise from physical considerations,
this article is concerned with the stabilization of
a linear wave equation subject to localized damping and a dynamical
boundary condition. This class of problems originates in
the classical treatment of the longitudinal and torsional vibrations of
rods. To illustrate
the physical relevance of this model, two representative settings are described below.
 
\begin{enumerate}
\item \textbf{Longitudinal impact of a rod.}
Consider a cylindrical rod occupying the interval $[0,L]$, fixed at the
end $x=0$, whose free end $x=L$ is struck at time $t=0$ by a point mass
$m$ moving with velocity $v$ along the axis of the rod. Let $u=u(x,t)$
denote the resulting longitudinal displacement of the rod. The governing
equation is the one-dimensional wave equation
$$
u_{tt} = c^2 u_{xx},
$$
where $c$ is the wave speed. The displacement satisfies the fixed-end condition $u(0,t)=0$, together with the dynamical boundary condition at
the impacted end,
$$
c^2 u_x(L,t) = -mL u_{tt}(L,t),
$$
which expresses Newton's second law for the impacting mass under the
reaction force transmitted by the rod. For a detailed derivation, see Ko\v{s}lyakov, Smirnov, and Gliner~\cite[p.~64]{Koshlyakov} and Timoshenko~\cite[p.~398]{Timoshenko}.
 
\item \textbf{Torsional vibration of a rod.}
A second example of the same structure occurs in the model for small
torsional deformations of a cylindrical rod fixed at the end
$x=0$. Let $\theta=\theta(x,t)$
denote the angle of torsion of the rod, $0\le x\le L$. The torsional displacement obeys
$$
\theta_{tt} = c^2 \theta_{xx},
$$
with the boundary condition $\theta(0,t)=0$ at the fixed end and the
dynamical boundary condition
$$
C^2 \theta_x(L,t) = -\theta_{tt}(L,t), \qquad 0<t<T,
$$
at the free end. Here $C$ represents the resistance to torsion, and $c$ is to be understood as in the
previous example. A detailed derivation can be found in Ko\v{s}lyakov, Smirnov, and Gliner~\cite[pp.~147--154]{Koshlyakov}.
\end{enumerate}

These two examples, although seemingly different in their physical
nature, are governed by the same abstract structure: a linear wave
equation on a one-dimensional domain, with a Dirichlet condition at
one end and a second-order dynamical boundary condition at the other.
The present article studies this structure in a broader setting:
 \begin{itemize}
     \item The wave equation is posed on a smooth multi-dimensional compact Riemannian manifold with
boundary. The boundary has at least two connected components: one where a Dirichlet boundary condition is imposed and another where a second-order dynamic boundary condition is used. 
  \item   A damping term is introduced either in the domain or at the boundary. If nonvanishing, this additional term yields decay of the energy.
 \end{itemize} 
In one dimension, a model that includes a damping term is derived in \cite[page 246]{rao2019vibration}. Higher-dimensional models with dynamical boundary conditions are derived for physical energy perspectives in \cite{Goldstein:06}.

 Before proceeding, the meaning of a \textit{dynamic boundary condition} is clarified. Unlike a static boundary condition, which imposes a constraint on the solution and its spatial derivatives at the boundary, a dynamic boundary condition prescribes an evolution equation for the boundary trace of the solution. In the present article, the undamped boundary dynamics are governed by
\begin{equation}\label{BC without Laplacian}
\d_\nu u +\d_t^2 u =0.
\end{equation}
Another type of dynamic boundary conditions that can be found in the literature is
\begin{equation}\label{BC with Laplacian}
\d_\nu u +\d_t^2u-\Delta_{\Gamma_1}u =0,
\end{equation}
where $\Delta_{\Gamma_1}$ denotes the Laplace-Beltrami operator of the connected component of the boundary $\Gamma_1$. There appears to be no universally accepted terminology distinguishing these two types of dynamic boundary conditions. Although these two types of boundary conditions reduce to the same one in the one-dimensional case (the tangential Laplacian $\Delta_{\Gamma_1}$ does not make sense anymore), they can lead to markedly different behaviors in higher dimensions. For instance, the stability results in case when damping occurs are different (see Section~\ref{sec:main_results}). From the controllability and stability perspective, boundary conditions \eqref{BC with Laplacian} have been extensively studied, but very few is known in the case of boundary conditions \eqref{BC without Laplacian} (see Section~\ref{Background} for such issue).

\subsection{Setting and functional framework}
Suppose $(\M,g)$ is a smooth $d$-dimensional Riemannian
manifold, and suppose $\d\M
\neq \emptyset$. In local coordinates, the Laplace-Beltrami operator $\Delta_g = \div_g \nablag$ is
\begin{align*}
\Delta_g f =(\det g)^{-1/2} 
  \sum_{1\leq i, j\leq d} \d_{x_i}\big( (\det g)^{1/2} g^{i j} \d_{x_j} f \big),
\end{align*}
where $g^* =(g^{ij})$ is the inverse of $g=(g_{ij})$.
Further assume that $\d\M$ has at least two connected
components and suppose that $\d\M = \Gamma_0\cup \Gamma_1$, where
$\Gamma_0$ and $\Gamma_1$ are disjoint and are nonempty unions of
connected components of $\d\M$. 

Suppose $\alpha \in L^\infty(\d\M)$, with $\alpha(x) \geq 0$ a.e. and set $\Ellop = - \Delta_g$, the problem studied in this article reads as follows
\begin{align}
  \begin{cases}
    \label{wp-1}
    \d_{t}^{2}u + \Ellop u=0
    \ \  \dans \ \ \M \times (0,+\infty), \\
    \dnu  u+\d_{t}^{2}u+\alpha \d_t u=0
    \ \ \sur \ \ \Gamma_1 \times (0,+\infty)  \ \ \et \ \
    u = 0 \ \ \sur \ \Gamma_0 \times (0,+\infty),\\
    u\ic=u_0, \quad \d_t u\ic=u_1
    \ \ \dans  \ \  \M.
  \end{cases}
\end{align}
At a point of $\Gamma_1$, the normal derivative $\dnu$ is 
understood in the sense of the metric; that is,
$\dnu u = g (\nablag u, \nu)$,
where $\nu$ is the unique unitary outward pointing vector in the sense of the metric.
In some places, we write $\gamma_0(u)$ for the Dirichlet trace on $\Gamma_0$, that is, $u_{|\Gamma_0}$. 

Denote by $\mu_g$ the Riemannian density on $\M$. For functions on $\M$, the
$L^2$-inner product and $L^2$-norm are  
\begin{align*}
  \inp{f}{f'}_{L^2(\M)} = \int_\M f \ovl{f'}\,  \mu_g
  \ \ \et \ \ 
  \Norm{f}{L^2(\M)} = \inp{f}{f}_{L^2(\M)}^{1/2}.
\end{align*}
For vector fields on $\M$, the
$L^2$-inner product and $L^2$-norm are  
\begin{align*}
  \inp{U}{U'}_{L^2 V(\M)} = \int_\M g(U, \ovl{U'})\,  \mu_g
  \ \ \et \ \ 
  \Norm{U}{L^2V(\M)} = \inp{U}{U}_{L^2 V(\M)}^{1/2}.
\end{align*}
Denote by $g_\d$ the induced metric on $\d\M$ and by $\mu_\d$ the
associated Riemannian density.  For functions on $\d\M$, and $\Gamma \subset \d\M$ an open set, the
$L^2$-inner product and $L^2$-norm on $\Gamma$ are
\begin{align*}
  \inp{f}{f'}_{L^2(\Gamma)} = \int_{\Gamma} f \ovl{f'}\, \mu_{\d}
  \ \ \et \ \ 
  \norm{f}{L^2(\Gamma)} = \inp{f}{f}_{L^2(\Gamma)}^{1/2}.
\end{align*} 

Introduce the following energy for a solution of \eqref{wp-1}
\begin{align}
  \label{eq:def_energy}
  \energy (u) = \frac12
  \big( \Norm{\nablag u}{L^2 V(\M)}^2 + \Norm{\d_t u}{L^2(\M)}^2
  + \norm{\d_t u}{L^2(\Gamma_1)}^2
  \big) .
\end{align}
It is a function of time. Consider the space
\begin{align}
  \label{eq:Hilbert_space_H}
  \H = H_{\Gamma_0}^{1}(\M) \times L^2(\M) \times  L^{2}(\Gamma_1),
\end{align}
where $H_{\Gamma_0}^{1}(\M)$ denotes the $H^1$ functions on $\M$ with vanishing trace on $\Gamma_0$. The space 
$\H$ is endowed with the inner product
\begin{align*}
  \inp{U}{V}_{\H} = \inp{\nablag u_{1}}{\nablag v_{1}}_{L^2 V(\M)}
  + \inp{u_2}{v_2}_{L^2(\M)}
  + \inp{u_{3}}{v_{3}}_{L^2(\Gamma_1)},
\end{align*}
for $U = (u_{1}, u_{2}, u_{3})$ and $V = (v_{1}, v_{2}, v_{3})$,
yielding a Hilbert structure on $\H$. Problem \eqref{wp-1} can be recast as 
\begin{align}
  \label{wp-3}
    \frac{d W}{d t}(t) =  \SgGen W(t), \qquad 
    W_{|t=0}= W_0= (u_{0}, v_{0}, {v_0}_{|\Gamma_1}),
\end{align}
 where $W=(u,\d_t u,\d_t u_{|\Gamma_1})$ and $\SgGen: \H \rightarrow \H$ is the unbounded operator given by
\begin{align}
    \label{wp-4}
    \SgGen=\begin{pmatrix} 0 & \id & 0 \\
    -\Ellop & 0 & 0 \\ -\dnu ( \cdot)_{|\Gamma_1} & -\alpha(x)(\cdot)_{|\Gamma_1} & 0	
    \end{pmatrix},
\end{align}
with domain
\begin{align}\label{wp-6}
	D(\SgGen)  =\Big\{ (u,v,w) \in \H; \ 
		v \in H_{\Gamma_0}^{1}(\M), \ 
		\Ellop u \in L^{2}(\M), \
		w= v_{|\Gamma_1}, \ 
		\dnu u_{|\Gamma_1} \in L^2(\Gamma_1)
	\Big\}.
\end{align}
For a solution $u$ to \eqref{wp-1} and the associated function $W(t)$,
define the energy of $W$ as the energy of $u$, that is, 
\begin{align*}
  \energy (W) = \energy (u) = \frac12 \Norm{W}{\H}^2.
\end{align*}

One says that  $W$ is a mild solution of \eqref{wp-3} if
\begin{align*}
W\in \Con([0,\infty);\mathcal{H}) \qquad \text{and} \qquad W(t)=e^{t \SgGen}W_0,
\end{align*}
where $e^{t \SgGen}$ is the semigroup generated by $A$.  This mild solution $W$ is called a strong
solution of \eqref{wp-3} if, in addition, 
\begin{align*}
W \in \Con([0,\infty); D(A)) \cap \Con^1([0,\infty); \mathcal{H}),    
\end{align*}
and
\begin{align*}
\frac{dW}{dt}(t) = AW(t)\quad \text{for all} \quad t \in [0,\infty), \qquad W(0) = W_0.
\end{align*}

For $k \in \N$, we define inductively $D\left(\SgGen^{k+1}\right)=
\{U \in D(\SgGen^k): \SgGen U \in D(\SgGen^k)\}$.

\begin{theorem}[Global well-posedness]\label{wp-thm-1}
  Suppose $W_0\in \H$. Then, System~\eqref{wp-3} has a unique
  mild solution $W \in \Con^0([0,\infty); \H)$.  Moreover, its energy
    $\energy (W)$ decays with time. If, in addition $U_0 \in D(\SgGen^k)$, then the mild solution is also the unique strong solution in $\bigcap_{j=0}^k \Con^j\big([0,\infty); D(\SgGen^{k-j})\big)$.
\end{theorem}

See Appendix~\ref{proof:wp-thm-1} for a proof.


\subsection{Main results}\label{sec:main_results}

Our main theorem states that there is at least a logarithmic decay of
the energy for the solutions of the system \eqref{wp-1}.
\begin{theorem}
  \label{thm:stabilization}
  Suppose $\Gamma$ is a nonempty open subset of $\Gamma_1$ and
  $\alpha \in \Con^{\infty}(\Gamma_1)$ is such that $\alpha>0$ on
  $\Gamma$. Suppose $k \in \N^*$. Then, there exists $C>0$ such
  that
  \begin{align}\label{ourdecay}
    \energy(W)(t) = \frac12 \Norm{W(t)}{\H}^2
    \leq \frac{C}{\log (2+t)^k} \Norm{\SgGen^k U^0}{\H}^2,
    \qquad \si \  W_0 \in D(\SgGen^k).
\end{align}
\end{theorem}
By the result of Batty–Duyckaerts~\cite{BattyDuyckaerts2008}, this
stability property is {\em equivalent} to the following resolvent
estimate for the operator $(\SgGen, D(\SgGen))$.
\begin{proposition}[resolvent estimate]
  \label{prop:resolvent_estimate}
Under the assumptions of Theorem~\ref{thm:stabilization}, the operator $i\sigma \id_{\H}-\SgGen$ is
invertible on $\H$ for all $\sigma \in \R$ and there exist $K>0$ and
$\sigma_0>0$ such that
\begin{align}
  \label{eq:resolvent_est}
  \Norm{(i \sigma \id_{\H}-\SgGen)^{-1}}{\mathscr{L}(\H, \H)}
  \leq K e^{K \sigma^2}, \quad \sigma \in \R, \ |\sigma| \geq \sigma_0 .
\end{align}
\end{proposition}
The resolvent estimate of Proposition~\ref{prop:resolvent_estimate} is
obtained by means of a Carleman estimate for a second-order elliptic
operator in which the spectral parameters appear as a potential and
also in the boundary condition inherited from \eqref{wp-1}. 

By a density argument, from Theorem~\ref{thm:stabilization} we obtain the following result.
\begin{corollary}
The energy of any mild solution to the boundary damped wave equation
\eqref{wp-3} goes to $0$ as $t \rightarrow+\infty$.
\end{corollary}

For the wave equation with the boundary condition of the form $\dnu
u+\alpha \d_t u=0$, a similar decay is known, yet with $\log(2+t)^{k}$
in \eqref{ourdecay} replaced by $\log(2+t)^{2k}$, that is a faster
decay; see \cite{Robbiano1997}, where optimality is also
proven. Under a dynamic boundary
  condition of type \eqref{BC with Laplacian}, \cite{Buffe2017} also
  proved a faster decay rate with $\log(2+t)^{2k}$. A natural
question is then the optimality of the result of
Theorem~\ref{thm:stabilization}. This is addressed in the next
proposition, via the optimality of the resolvent
estimate.
\begin{proposition}[optimality of the resolvent estimate]
  \label{prop:optimality}
  The factor $\sigma^2$ in the resolvent estimate~\eqref{eq:resolvent_est} is optimal, that is, it cannot be replaced by $\sigma^\eta$ for any $0<\eta<2$. More precisely, for any $d \geq 2$, there exist a Riemannian manifold $(\M,g)$ with boundary, a constant $C > 0$, a damping coefficient $\alpha$, and a sequence 
  $(\sigma_n)_{n\in\N}$ of real numbers with $\sigma_n \to +\infty$, such that
  \begin{align*}
    \Norm{(i \sigma_n \id_{\H}-\SgGen)^{-1}}{\H}
    \geq C e^{C \sigma_n^2}, \quad n \in \N.
  \end{align*}
\end{proposition}
This optimality result emphasizes the fact that the dynamic boundary
condition of System \eqref{wp-1} is fundamentally different from the Neumann condition
\cite{Robbiano1997} or a dynamic boundary condition of the type given in 
\eqref{BC with Laplacian}.  This is furthermore discussed in
Section~\ref{sec:optimiality_Carleman} in connection with the Carleman
estimates that are derived in Section~\ref{sec:Carleman_estimate} and
that lead to the proof of the resolvent estimate of
Proposition~\ref{prop:resolvent_estimate}.
  
Results similar to those of Theorem~\ref{thm:stabilization} and
Propositions~\ref{prop:resolvent_estimate} and \ref{prop:optimality}
are also obtained in the case of interior damping; see
Section~\ref{sec:interior}.

\subsection{Background} \label{Background}
The control and the stabilization of PDEs have become central
questions in their study, in particular in the case of the wave
equation.  The multiplier method played an early and crucial role; see
\cite{Komornik1994, Lionsvol1, Lionsvol2}. This technique allows one
to derive energy estimates and to establish observability inequalities,
and, as a consequence, many controllability results by means of the
duality method, which was popularized by Lions under the name HUM,
standing for Hilbert Uniqueness Method \cite{Ho1986, Lionsvol1}.

Microlocal analysis came as another powerful approach allowing one to
describe how waves interact with the geometry of the domain, leading
to refined results in control and inverse problems.  In this context,
we mention the seminal work of Bardos, Lebeau, and Rauch
\cite{BLR1992} which establishes sharp geometric sufficient
conditions for these properties to hold when the control or
observation is applied from the boundary of a domain. Their central
contribution was the introduction of the Geometric Control Condition,
which reveals a connection between the observation region
and the generalized geodesics of the wave operator.

Carleman estimates have also emerged as an important tool in the
analysis of partial differential equations. These weighted inequalities are
particularly effective in proving unique continuation properties and
have become an efficient method in control theory and inverse
problems. Closely related to these developments are resolvent
estimates for the wave equation. They provide important insights into the spectral behavior
of the wave semigroup and play a key role in the understanding of stabilization; see \cite{Lebeau1994, Robbiano1997, Robbiano1995}. The present contribution lies in this line of development.

In recent years, evolution equations with dynamic boundary conditions have attracted increasing attention due to their relevance in several applied contexts. These boundary conditions have been the focus of significant research in the context of control and stabilization problems for partial differential equations. Several works address the stabilization of the wave equation under such boundary conditions.

Regarding the available literature on dynamic boundary conditions in problem \eqref{wp-1}, all results concern damping that acts on the whole domain. We first mention Araújo, Milla Miranda, and Medeiros
\cite{Milla2009MC}, who investigate the well-posedness and long-time
behavior of a damped wave equation $\d_t^2u-\Delta u +a(x)\d_t u=0$ posed
on a bounded domain $\Omega \subset \R^n$ with a dynamic
boundary condition of the form $\d_\nu u + \beta(x)\d_t^2u +
\alpha(x)\d_t u = 0$ on the connected component $\Gamma_1$ of the boundary. Defining $H_{\Gamma_0}^{1}(\Omega)=\{u\in H^1(\Omega): u_{|\Gamma_0}=0\}$ endowed with the norm $\Norm{u}{H_{\Gamma_0}^{1}(\Omega)}^2=\int_{\Omega}|\nabla u|^2dx$, the authors
prove the existence and uniqueness of strong solutions. Under the
additional assumptions $a(x) \geq
a_0 > 0$, $\alpha(x) \geq \alpha_0 > 0$,  and $\beta \neq 0$, it is shown via a perturbed energy argument that the energy, defined as $2E(t) = \|\d_t u(t)\|^2_{L^2(\Omega)} +
\|u(t)\|^2_{H_{\Gamma_0}^{1}(\Omega)} + \norm{\beta^{1/2}\d_t u(t)}{L^2(\Gamma_1)}^2$,  decays exponentially at the rate $e^{-\frac{2}{3}\eta t}$ where $\eta>0$. For further results on well-posedness and stabilization obtained by the same approach, see also \cite{Milla2024MMAS, Milla2022MC}. 

 Milla Miranda, Medeiros, and Lourêdo \cite{Milla2022EJDE} investigate the
existence, uniqueness, and exponential decay of solutions for a
quasilinear hyperbolic problem posed on a bounded domain $\Omega
\subset \R^n$, motivated by the nonlinear longitudinal
vibrations of a bar with a concentrated mass attached at one end. This
model is deduced from Newton's laws and the nonlinear Hooke's law
$\tau(x,t) = \sigma(u_x(x,t))$; see \cite{Milla2024MMAS}. The problem reads 
$\d_t^2 u - \sum_{i=1}^{n} \d_{x_i}\big[\sigma_i(\d_{x_i} u)+\d_{x_i}(\d_t u)\big]=0$ 
in $\Omega$, with $u=0$ on $\Gamma_0$, 
and featuring a dynamic boundary
condition of the form $\sum_{i=1}^n \bigl[\sigma_i(\d_{x_i} u) +
  \d_{x_i} \d_t u\bigr]\nu_i + \d_t^2u = 0$ on $\Gamma_1$. The authors point out that
global existence for the analogous problem with zero Dirichlet
boundary conditions and without internal damping remains open, which
motivates the introduction of internal dissipation. Under the growth condition $s^2 \leq b_i\,\widehat{\sigma}_i(s)$ for each $i=1,\dots,n$, where $b_i>0$ are constants and $\widehat{\sigma}_i(s) = \int_0^s \sigma_i(\tau)\,d\tau$, the authors prove that the energy
$2E(t) = \Norm{\d_t u(t)}{L^2(\Omega)}^2 + \sum_{i=1}^n \int_\Omega \widehat{\sigma}_i(\d_{x_i}u(t))dx + \norm{\d_t u(t)}{L^2(\Gamma_1)}^{2}$ 
decays exponentially at the rate $e^{-\frac{2}{3}\eta t}$ where $\eta>0$, via Lyapunov's approach. 

 Li, Liang, and Xiao \cite{Li2021} study the asymptotic behavior of the
linear wave equation with frictional damping appearing only on the
dynamical boundary, i.e., $u_{tt} + \dnu u + u_t = 0$ on
$\Gamma_1$, without any interior damping. By analyzing the spectrum of the semigroup generator $\SgGen$, they show that
$\sigma(\SgGen) \cap i\R = \emptyset$ and that $\SgGen^{-1}$ is
compact, which implies strong stability of the energy via the
Arendt–Batty–Lyubich–Vũ theorem. To obtain quantitative decay rates,
the authors use multiplier methods to derive resolvent estimates. Then, by the Borichev–Tomilov theorem \cite{Borichev2010}, they establish polynomial decay of the energy, which they prove to be optimal in the 1D case \cite{Li2021}. 

In a related direction, Milla Miranda and Lour\^edo \cite{Milla2023EECT} studied the exact controllability of the one-dimensional wave equation under the dynamic boundary condition $\d_t^2 u(L,t) + \d_x u(L,t) = 0$ through the HUM method, thus complementing the stabilization results discussed above for this class of boundary conditions.

In the context of damped wave equations with dynamic boundary conditions of type \eqref{BC with Laplacian}, logarithmic decay and resolvent estimates are derived in \cite{Buffe2017} by means of Carleman estimates. In addition, one may also refer to \cite{Cavalcanti2007, Cavalcanti2012, Parada2025} for further related results.  In particular, optimal polynomial decay has recently been proven in \cite{Parada2025} even under a GCC assumption.

When considering purely reactive dynamic boundary conditions such as
those considered in this article, previous studies only focused on
\emph{full damping} settings, where the dissipation acts either
throughout the whole domain in the case of internal damping
\cite{Milla2009MC, Milla2024MMAS, Milla2022EJDE}, or on the entire
connected component $\Gamma_1$ in the case of boundary damping
\cite{Li2021, Li2017}. The situation where dissipation is confined to
a small region, whether in the interior or on a portion of the
boundary, and without any geometric restriction, has not been
addressed before. This is the gap that the present work seeks to
fill. In fact, we prove that the presence of the inertial term
$\d_t^2 u$ together with a localized damping term alters, in an
essential way, the qualitative behavior of the solutions, as already
pointed out in the discussion preceding
Proposition~\ref{prop:optimality}.

\subsection{Outline}

The remainder of the present paper is organized as follows. Section~\ref{sec:psido} 
introduces some material on pseudo-differential operators with large parameters and adapted Sobolev norms. 
In Section~\ref{sec:Carleman_estimate}, we
establish a local Carleman estimate at the boundary. After introducing the local setting in
normal geodesic coordinates, a conjugation is performed using an
exponential weight function, and the resulting operators  are analyzed. The proof proceeds through a microlocal analysis under the \LS condition, distinguishing two cases depending on
the location of the roots of the principal symbol. The key estimates
are obtained via a microlocal parametrix, culminating in the proof of the Carleman estimate.

In Section~\ref{sec:spectral_resolvent_estimate}, the resolvent estimate for the generator $\SgGen$
of the semigroup associated with problem
\eqref{wp-1} is derived. Using the Carleman estimate established in
Section~\ref{sec:Carleman_estimate}, it is proved that the resolvent grows at most exponentially in
$\sigma^2$, specifically $\|(i\sigma \id_{\H} -
\SgGen)^{-1}\|_{\mathcal{L}(\H)} \leq K e^{K\sigma^2}$ for $\sigma \in \R$. This estimate implies a logarithmic decay of the wave semigroup, yielding our main stabilization
result in Theorem~\ref{thm:stabilization}. Section \ref{sec:interior} presents analogous results for the case of interior
damping. Finally, in Section \ref{sec:optimality}, the optimality of the resolvent estimate
is established in well-chosen geometries by constructing quasi-modes that concentrate
away from the damping region, showing that the exponent $2$ in the
resolvent bound cannot be improved.

\subsection{Notation}
Write
\begin{align*}
  \Rdp:=\{x \in \R^d; x_d>0\} \quad \where\  x=(x', x_d) \  \avec\  x' \in \R^{d-1}, x_d \in \R.
\end{align*}
Use accordingly $\xi=(\xi', \xi_d) \in \R^{d-1} \times \R$ and
$D_d=-i \d_d$ and $D'=-i \d'$, with $\d'=(\d_{x_1}, \dots,
\d_{x_{d-1}})$.


For an open set $U$ of $\R^d$, set $U_{+}=U \cap \Rdp$ and
\begin{align*}
\Cbarc(U_{+})=\big\{ u=v_{|\Rdp};\  v \in \Cinfc(\R^d)  \ \et \  \supp v \subset U\big\}.
\end{align*}
Set also $\ovl{\S}(\Rdp)=\left\{u_{| \Rdp} ; u \in
\S\left(\R^d\right)\right\}$ with $\S(\R^d)$ the
usual Schwartz space in $\R^d$.

Recall that the Poisson bracket of two smooth functions is given by
\begin{align*}
\{f, g\}= \sum_{j=1}^d (\d_{\xi_j} f \d_{x_j} g-\d_{x_j} f \d_{\xi_j} g).
\end{align*}

The notation $a \lesssim b$ is used for $a \leq C b$ and $a \gtrsim b$ for $a \geq C b$, with a constant $C>0$ that may change from one line to another. Moreover, $a \eqsim b$ denotes $a \lesssim b \lesssim a$.

\subsection{AI disclaimer}
The mathematical analysis presented in this article was carried without any use of AI tools.

\section{Pseudo-differential operators with large parameters}
\label{sec:psido}

We recall some notions on semi-classical pseudo-differential operators
with a large parameter $\tau \geq 1$. Set $\lambda_\tau^2=\tau^2+|\xi|^2$.
For $m \in \R$, denote by
$S_\tau^m$ the space of smooth functions $a(x, \xi, \tau)$ defined for
$(x,\xi)\in \R^d \times \R^d$ and for $\tau \geq 1$ such that for all
multi-indices $\alpha, \beta \in \N^d$, for some 
$C_{\alpha, \beta}>0$ one has
\begin{align*}
  |\d_x^\alpha \d_{\xi}^\beta a(x, \xi, \tau)|
  \leq C_{\alpha, \beta} \lambda_\tau^{m-|\beta|}, \qquad (x, \xi, \tau) \in \R^d \times \R^d \times[1, \infty).
\end{align*}
If $a \in S_\tau^m$, define the associated pseudo-differential operator
\begin{align*}
\Op(a) u(x) :=\frac{1}{(2 \pi)^d} \int_{\R^d} e^{i x \cdot \xi} a(x, \xi, \tau) \hat{u}(\xi) d \xi, \quad u \in \S(\R^d).
\end{align*}
One says that $a$ is the symbol of $\Op(a)$ and that $\Op(a)$ is of
order $m$.  Denote $\Psisc^m$ the set of such operators.  Denote
$\mathscr{D}_\tau^m$ the space of semi-classical differential
operators, i.e., the case when the symbol $a(x, \xi, \tau)$ is a
polynomial function of order $m$ in $(\xi,\tau)$.
If $a\in S_\tau^m$, then the principal symbol of $\Op(a)$ is (a representative of) the class of $a$ in $S_\tau^m / S_\tau^{m-1}$.

\subsection{Tangential operators}
Introduce pseudo-differential operators that act in the tangential
direction $x'$ with $x_d$ behaving as a parameter as follows.  Set $\lsct^{2} = \tau^2+|\xi'|^2$. Denote
$\Ssct^{m}$, the set of smooth functions $b(x, \xi', \tau)$ defined
for $(x,\xi')\in \R^d \times \R^{d-1}$ and for $\tau \geq 1$ such that  for all
multi-indices $\alpha \in \N^d, \beta \in \N^{d-1}$, for some 
 $C_{\alpha, \beta}>0$ one has 
\begin{align}
  \label{eq:tang_symbols}
  |\d_x^\alpha \d_{\xi'}^\beta b (x, \xi', \tau)|
  \leq C_{\alpha, \beta} \lsct^{m-|\beta|},  \qquad  (x, \xi', \tau) \in \R^d \times \R^{d-1} \times[1, \infty).
\end{align}
If $b \in \Ssct^m$, define the associated pseudo-differential operator
\begin{align*}
  \opt(b) u(x):=\frac{1}{(2 \pi)^{d-1}} \int_{\R^{d-1}} e^{i x' \cdot \xi'} b(x, \xi', \tau) \hat{u}(\xi', x_d) d \xi',
  \quad u \in \ovl{\S}(\R^d_+).
\end{align*}
Denote $\Psisct^{m}$  the set of such tangential pseudo-differential operators of order $m$. 
Here, a subset of $\R^d \times\R^{d-1}\times [1,\infty)$ is said to be conic if it is conic in the $(\tau,\xi)$ variables.
If $\U$ is such a conic open set, one says that $a \in \Ssct^{m}$ microlocally in $\U$ if
  \eqref{eq:tang_symbols} holds for $(x,\xi',\tau)\in \U$.

Suppose $m \in \N$, $r \in \R$, and $a$ is of the form
\begin{align*}
  a(x, \xi, \tau)=\sum_{j=0}^m a_j(x, \xi', \tau) \xi_d^j, \quad \avec \  a_j \in \Ssct^{m+r-j}.
\end{align*}
Define $\Op(a) = \sum_{j=0}^{m} \Opt(a_j) D_{x_d}^j$.
One says that $a \in \Ssc^{m, r}$ and $\Op(a) \in \Psisc^{m, r}$.

\begin{proposition}[calculus results]
  \label{prop:calculus}
  Suppose  $a \in \Ssct^m$ and $b \in \Ssct^{m'}$.
  Then,
  \begin{enumerate}
  \item $\Opt(a)^* \in \Psisc^m$ with $\bar{a}$ for principal symbol;
  \item $ab \in \Ssct^{m+m'}$ and $\Opt(a)\Opt(b) \in \Psisc^{m+m'}$ with $ab$  for principal symbol;
  \item $[\Opt(a), \Opt(b)] \in  \Psisc^{m+m'-1}$ with $\frac{1}{i} \{a,b \} \in \Ssct^{m+m'-1}$ for principal symbol. 
  \end{enumerate}
\end{proposition}
We refer to \cite[Chapter~2]{JLR1} for a proof.

\subsection{Sobolev norms and bounds}
\label{sec:Sobolev_norms_bounds}
We use the notation $\Norm{.}{}$ for functions defined in the interior
of the domain and $\norm{.}{}$ for functions defined on the
boundary. For example, set 
\begin{align*}
\Norm{u}{+}=\Norm{u}{L^2(\Rdp)} \quad  \et \quad \inp{u}{\tilde{u}}_{+}=\inp{u}{\tilde{u}}_{L^2(\Rdp)},
\end{align*}
for functions defined in $\Rdp$ and
\begin{align*}
  \norm{w}{\d}=\norm{w}{L^2(\R^{d-1})}  \quad  \et \quad \inp{w}{\tilde{w}}_{\d} = \inp{w}{\tilde{w}}_{L^2(\R^{d-1})},
\end{align*}
for functions defined on $\{x_d=0\}$, such as traces.

For $r\in \R$, set $\Lsct^{r}=\Opt(\lsct^{r})$. For $m \in \N$ set
\begin{align*}
  \Normsc{v}{m,r} =\sum_{j=0}^m\Norm{\Lsct^{m-j+r} D_d^j v}{+}
  \quad \et \quad  \Normsc{v}{m} = \Normsc{v}{m,0}
  =\tau^{m} \Norm{v}{+} + \Norm{v}{H^m(\Rdp)}.
\end{align*}
For $s \in \R$, set 
\begin{align*}
  \normsc{w}{s} = \norm{\Lsct^s w}{\d}.
\end{align*}
If $s \geq 0$, then $\normsc{w}{s} = \tau^{s} \Norm{w}{\d}
  + \Norm{w}{H^s(\R^{d-1})}$.
For $m \in \N$, and $s \in \R$, set 
\begin{align*}
  \normsc{\trace(v)}{m,s}
  = \sum_{j=0}^m
  \normsc{D_{d}^j v\br}{m+s-j}.
\end{align*}
\begin{remark}\label{remark: tau parameter suff large ch6}
Observe that, for some $C>0$, $\Normsc{v}{m,r} \leq C
\tau^{-\ell}\Normsc{v}{m,r+\ell}$ for $u \in \ovl{\S}(\Rdp)$, $m \in
\N$, $s \in \R$, and $\ell \geq 0$. This implies that $\Normsc{v}{m,r}
\ll\Normsc{v}{m,r+\ell}$ for $\tau$ large.
\end{remark}

For $r\in \R$, set also $\Lsc^{r}=\Op(\lsc^{r})$.
\begin{proposition}[operators bounds]
  \label{prop:operator_bounds}
  Suppose $r,s \in \R$.

  \begin{enumerate}
    \item Suppose $a \in \Ssc^s$. Then, there exists $C>0$ such that 
      $\Norm{\Lsc^{r}\Op(a) u}{L^2(\R^d)}\leq C \Norm{\Lambda^{r+s}u}{L^2(\R^d)}$ for any $u \in \S(\R^d)$.
      
    \item Suppose $b \in \Ssct^s$. Then, there exists $C>0$ such that 
      $\Normsc{\Opt(a) u}{r} 
        \leq C \Normsc{u}{r+s}$ for any $u \in \ovl{\S}(\Rdp)$.
     \item Suppose $m,\ell\in \N$ and $c\in \Ssct^{m,s}$.
      Then, there exists $C>0$ such that 
      \begin{align*}
        \Normsc{\Opt(c) u}{\ell,r} 
        \leq C \Normsc{u}{m+\ell, r+s}, \qquad u \in \ovl{\S}(\Rdp).
      \end{align*}
  \end{enumerate}
\end{proposition}
We refer to \cite[Section~6.2]{JLR2} for a proof. 

The following trace inequalities are adapted to the norms we introduced. 
\begin{proposition}[trace inequalities]
  \label{prop:trace_inequality}
  Suppose $m\in \N$ and $s>0$.
  There exists $C>0$ such that
  \begin{align*}
    \normsc{\trace(u)}{m,s} \leq C \Normsc{u}{m+1, s-1 / 2}, \quad u \in \ovl{\S}(\Rdp).
  \end{align*}
\end{proposition}

\subsection{Symbols and operators with an additional parameter}

In this article, we also use operators and symbols that depend on an
additional parameter $\sigma \in \R$, which stands for the spectral
parameter in our analysis.

Set $\lscbis{\sigma}^2 = \tau^2+\sigma^2+|\xi|^2$ and $\lsctbis{\sigma}^2 = \tau^2+\sigma^2+|\xi'|^2$.
We are interested by symbols $a(x, \xi, \tau, \sigma)$ and tangential symbols $b(x, \xi', \tau, \sigma)$  that fulfill 
\begin{align*}
  |\d_x^\alpha \d_{\xi}^\beta a(x, \xi, \tau, \sigma)| \leq C_{\alpha, \beta}\lscbis{\sigma}^{(m-|\beta|) / 2}
\ \ \et \ \
  |\d_x^\alpha \d_{\xi'}^\beta b(x, \xi', \tau, \sigma)| \leq C_{\alpha, \beta}\lsctbis{\sigma}^{(m-|\beta|) / 2}.
\end{align*}
Observe that if $\tau \gtrsim |\sigma|$ one has
\begin{align*}
  \lscbis{\sigma} \eqsim \lsc \ \ \et \ \
   \lsctbis{\sigma} \eqsim \lsct.
\end{align*}
Thus, as far as pseudo-differential calculus and Sobolev
bounds are concerned, it is as if $a \in \Ssc^m$ and $b\in \Ssct^m$. As one always assumes that $\tau
\gtrsim |\sigma|$, this property
shall be exploited throughout the article.

\section{Local Carleman estimates at the boundary}
\label{sec:Carleman_estimate}

We prove Carleman estimates for the operator $\Psig = \Ellop -
\sigma^2$ where $\sigma$ is the spectral parameter in the resolvent
estimate of Proposition~\ref{prop:resolvent_estimate}. The parameter $\sigma$ is thus assumed to be real with $|\sigma| \gtrsim 1$.

These estimates are first proven in the case of Neumann boundary conditions
and then adapted to the boundary conditions needed for the treatment of
the present operators. 

In agreement with the norms on $\Rdp$ given in Section~\ref{sec:Sobolev_norms_bounds}, define the following norms on $\M$ for $m \in \N$ and $s >0$
\begin{align*}
  \Normsc{v}{m} =
  \tau^{m} \Norm{v}{L^2(\M)}
  + \Norm{v}{H^m(\M)}
  \ \ \et \ \
  \normsc{v_{|\d\M}}{s} =
  \tau^{s} \Norm{v_{|\d\M}}{L^2(\d\M)}
  + \Norm{v_{|\d\M}}{H^s(\d\M)}.
\end{align*}
Then, for $0\leq r \leq m$, one has $\tau^{m-r} \Norm{v}{H^r(\M)}\lesssim \Normsc{v}{m}$.
 For $m \in \N$, and $s \in \R$ and a smooth function $f$ defined on $\d \M$, one defines
\begin{align*}
  \normsc{f\trace(v)}{m,s}
  = \sum_{j=0}^m
  \normsc{f\dnu^j v_{|\d\M}}{m+s-j}.
\end{align*}

The principal symbol of $\Ellop$ is $\ellop(x,\xi) =  g_x^*(\xi,\xi)$.
Suppose $x^0 \in \d\M$ and $\varphi$ is a smooth function defined in a
\nhd $V$ of $x^0$. Introduce the (principal) symbol $\psigphi (x,\xi,\tau) =
\ellop(x,\xi + i \tau d \varphi(x)) - \sigma^2$, written here in local coordinates. Recall the following definition.
\begin{definition}[\se]
  \label{def:sub_ellipticity}
  One says that $(\Psig, \varphi)$ has the \se
  property in $\ovl{V}$ if $|d \varphi|>0$ in $\ovl{V}$ and 
  \begin{align*}
    \forall (x,\xi) \in \ovl{V} \times \R^d, \ \forall \tau>0,
    \qquad \psigphi (x,\xi,\tau) =0 \ \ \imply \ \
    \frac{1}{2i} \{ \ovl{\psigphi}, \psigphi\}(x,\xi,\tau) >0.
  \end{align*}
\end{definition}
Examples of functions $\varphi$ that fulfill this property are recalled in Section~\ref{sec:local_coordinates}. 
\begin{theorem}
   \label{th:Carleman_boundary1}
   Suppose $\dnu \varphi <0$ in $\ovl{V} \cap \d\M$.  Suppose $C_0>0$. Suppose that
   $(P, \varphi)$ has the \se property in $\ovl{V}$.  Then, there
   exist $C>0$, $\tau_*>0$, and a \nhd $W$ of $x^0$, with $W\subset V$, such that
   \begin{align}
     \label{eq:Carleman_boundary1}
     \tau^{-1/2} \Normsc{e^{\tau \varphi} u}{2}
     + \normsc{e^{\tau \varphi} \trace(u)}{1,1/2}
     \leq C \big( \Norm{e^{\tau \varphi} \Psig u}{L^2(\M)}
     + \normsc{e^{\tau \varphi} \dnu u}{1/2}\big),
   \end{align}
for $|\sigma| \geq 1$, $\tau \geq \max(\tau_*, C_0|\sigma|)$, and $u \in
\Cinf(\M)$ with $\supp u \subset W$.
\end{theorem}
If no assumption is made on the sign of $\dnu \varphi$ at the
boundary, one has the following estimate.
\begin{theorem}
   \label{th:Carleman_boundary2}
   Suppose that
   $(P, \varphi)$ has the \se property in $\ovl{V}$. Suppose $C_0>0$.  Then, there
   exist $C>0$, $\tau_*>0$, and a \nhd $W$ of $x^0$, with $W\subset V$, such that
   \begin{align}
     \label{eq:Carleman_boundary2}
     \tau^{-1/2} \Normsc{e^{\tau \varphi} u}{2}   
     \leq C \big( \Norm{e^{\tau \varphi} \Psig u}{L^2(\M)}
     + \normsc{e^{\tau \varphi} \trace(u)}{1,1/2}\big),
   \end{align}
for $|\sigma| \geq 1$, $\tau \geq \max(\tau_*, C_0|\sigma|)$, and $u \in
\Cinf(\M)$ with $\supp u \subset W$.
\end{theorem}
Sections~\ref{sec:local_coordinates} to
\ref{sec:microlocal_estimates} are dedicated to the proof of
Theorems~\ref{th:Carleman_boundary1} and \ref{th:Carleman_boundary2}.

Adapting {\em mutatis mutandis} the proof of Theorem~8.24 in
\cite{JLR2}, one deduces the
following shifted estimate from Theorem~\ref{th:Carleman_boundary1}.
\begin{corollary}
   \label{cor:Carleman_boundary-shift1}
   Suppose $\dnu \varphi <0$ in $\ovl{V} \cap  \d\M$. Suppose $C_0>0$.  Suppose that
   $(P, \varphi)$ has the \se property in $\ovl{V}$.  Then, there
   exist $C>0$, $\tau_*>0$, and a \nhd $W$ of $x^0$, with $W\subset V$, such that
   \begin{align}
     \label{eq:Carleman_boundary-shift1}
     \tau^{1/2} \Normsc{e^{\tau \varphi} u}{1}
     + \tau^{1/2}  \normsc{e^{\tau \varphi} \trace(u)}{1,0}
     \leq C \big( \Norm{e^{\tau \varphi} \Psig u}{L^2(\M)}
     + \tau^{1/2} \norm{e^{\tau \varphi} \dnu u}{L^2(\d \M)}\big),
   \end{align}
for $|\sigma| \geq 1$, $\tau \geq \max(\tau_*, C_0|\sigma|)$, and $u \in
\Cinf(\M)$ with $\supp u \subset W$.
\end{corollary}
Similarly, adapting  the proof of Theorem~8.30 in
\cite{JLR2}, one deduces the
following shifted estimate from Theorem~\ref{th:Carleman_boundary2}.
\begin{corollary}
   \label{cor:Carleman_boundary-shift2}
   Suppose that
   $(P, \varphi)$ has the \se property in $\ovl{V}$.  Suppose $C_0>0$.  Then, there
   exist $C>0$, $\tau_*>0$, and a \nhd $W$ of $x^0$, with $W\subset V$, such that
   \begin{align}
     \label{eq:Carleman_boundary-shift2}
     \tau^{1/2} \Normsc{e^{\tau \varphi} u}{1}   
     \leq C \big( \Norm{e^{\tau \varphi} \Psig u}{L^2(\M)}
     + \tau^{3/2}\norm{e^{\tau \varphi} u_{|\d\M}}{L^2(\d\M)}
     + \tau^{1/2}\norm{e^{\tau \varphi} \dnu u_{|\d\M}}{L^2(\d\M)}
     \big),
   \end{align}
for $|\sigma| \geq 1$, $\tau \geq \max(\tau_*, C_0|\sigma|)$, and $u \in
\Cinf(\M)$ with $\supp u \subset W$.
\end{corollary}

We now provide similar statements, yet by replacing the Neumann trace
on the \rhs of the estimates by the trace operator \
\begin{equation}
    \label{eq:Bsig}
    \Bsig= \dnu - \sigma^2 + i \sigma \alpha.
\end{equation} In fact, as $\alpha$
is bounded, for $|\sigma| \geq 1$ one finds
\begin{align*}
  \tau^{1/2}\norm{e^{\tau \varphi} \dnu u_{|\d\M}}{L^2(\d\M)}
  \lesssim
  \tau^{1/2}\norm{e^{\tau \varphi} \Bsig u_{|\d\M}}{L^2(\d\M)}
  + \tau^{1/2}\sigma^2 \norm{e^{\tau \varphi} u_{|\d\M}}{L^2(\d\M)}
\end{align*}

The second term on the \rhs can be absorbed by the \lhs of
\eqref{eq:Carleman_boundary-shift1} if $\tau \sigma^{-2}$ is chosen \suff
large. From Corollary~\ref{cor:Carleman_boundary-shift1} one obtains
the following proposition.
\begin{proposition}
   \label{prop:Carleman_boundary-Bsigma-1}
   Suppose $\dnu \varphi <0$ in $\ovl{V} \cap \d\M$.  Suppose that
   $(P, \varphi)$ has the \se property in $\ovl{V}$.  Then, there
   exist $C>0$, $\tau_*>0$, and a \nhd $W$ of $x^0$, with $W\subset V$, such that
   \begin{align}
     \label{eq:Carleman_boundary-Bsigma-1}
     \tau^{1/2} \Normsc{e^{\tau \varphi} u}{1}
     + \tau^{1/2}  \normsc{e^{\tau \varphi} \trace(u)}{1,0}
     \leq C \big( \Norm{e^{\tau \varphi} \Psig u}{L^2(\M)}
     + \tau^{1/2}\norm{e^{\tau \varphi} \Bsig u_{|\d \M}}{L^2(\d\M)}\big),
   \end{align}
for $|\sigma| \geq 1$, $\tau \geq \tau_* \sigma^2$, and $u \in
\Cinf(\M)$ with $\supp u \subset W$.
\end{proposition}
Similarly, working with the r.h.s of \eqref{eq:Carleman_boundary-shift2} in Corollary~\ref{cor:Carleman_boundary-shift2}, one
obtains the following proposition.
\begin{proposition}
   \label{prop:Carleman_boundary-Bsigma-2}
   Suppose that
   $(P, \varphi)$ has the \se property in $\ovl{V}$.  Then, there
   exist $C>0$, $\tau_*>0$, and a \nhd $W$ of $x^0$, with $W\subset V$, such that
   \begin{align}
     \label{eq:Carleman_boundary-Bsigma-2}
     \tau^{1/2} \Normsc{e^{\tau \varphi} u}{1}   
     \leq C \big( \Norm{e^{\tau \varphi} \Psig u}{L^2(\M)}
     + \tau^{3/2}\norm{e^{\tau \varphi} u_{|\d\M}}{L^2(\d\M)}
     + \tau^{1/2}\norm{e^{\tau \varphi} \Bsig u_{|\d\M}}{L^2(\d\M)}
     \big),
   \end{align}
for $|\sigma| \geq 1$, $\tau \geq \tau_* \sigma^2$, and $u \in
\Cinf(\M)$ with $\supp u \subset W$.
\end{proposition}

\subsection{Optimality issues}
\label{sec:optimiality_Carleman}
Estimate \eqref{eq:Carleman_boundary-Bsigma-1} in
Proposition~\ref{prop:Carleman_boundary-Bsigma-1} is deduced from
Corollary~\ref{cor:Carleman_boundary-shift1}. Estimate
\eqref{eq:Carleman_boundary-Bsigma-2} in
Proposition~\ref{prop:Carleman_boundary-Bsigma-2} is deduced from
Corollary~\ref{cor:Carleman_boundary-shift2}. These two steps are the
only parts of the argumentation where the condition $\tau \gtrsim
|\sigma|$ has to be replaced by the {\em much stronger} condition
$\tau  \gtrsim \sigma^2$. This has a great impact on the resolvent
estimate of Proposition~\ref{prop:resolvent_estimate} through the
factor $\exp(C \sigma^2)$. The condition $\tau \gtrsim |\sigma|$ would
yield a factor $\exp(C |\sigma|)$, similar to what one finds in the
case of damping through the boundary operator $\dnu + \alpha \d_t$;
see \cite{Robbiano1997} and \cite[Chapter~10]{JLR2}.

  The question of the optimality of the condition $\tau \gtrsim
  \sigma^2$ in the Carleman estimates of
  Propositions~\ref{prop:Carleman_boundary-Bsigma-1}
  and~\ref{prop:Carleman_boundary-Bsigma-2} thus appears very natural.
  The answer is, in fact, {\em positive}: in
  Section~\ref{sec:optimality}, we prove that the factor $\exp(C
  \sigma^2)$ in the resolvent estimate cannot be improved. Hence, the
  condition $\tau \gtrsim \sigma^2$ cannot be improved. It is {\em not} a
  limitation of the method. The difference between
  Corollary~\ref{cor:Carleman_boundary-shift1} and
  Propositon~\ref{prop:Carleman_boundary-Bsigma-1} on the one hand,
  and the difference between
  Corollary~\ref{cor:Carleman_boundary-shift2} and
  Proposition~\ref{eq:Carleman_boundary-Bsigma-2} on the other hand,
  clearly show that the dynamic boundary condition of System \eqref{wp-1}
  is fundamentally different from the Neumann condition or the dynamic
  boundary condition of type \eqref{BC with Laplacian}.

\subsection{Local coordinates and operator conjugation}
\label{sec:local_coordinates}

Suppose $x^0 \in\d\M$ and $V$ an open set of $\M$, a \nhd of $x^0$,
where the \se property holds as in the statement of
Theorems~\ref{th:Carleman_boundary1} and
\ref{th:Carleman_boundary2}. Then, there exists a \nhd $W \subset V$
of $x^0$ where one can consider normal geodesic coordinates. In such
coordinates, one has $\d\M = \{ x_d=0\}$, $\M = \{ x_d >0\}$, and the
metric $g$ takes the form
\begin{align*}
  g = g' + d_{x_d} \otimes d_{x_d}
  \ \ \avec \ \
  g' = \sum_{1\leq i, j\leq d-1} g_{ij}'\,  d_{x_i} \otimes d_{x_j}.
\end{align*}
In these coordinates $\nu = - \d_{x_d}$.
Classically, we set $D = - i \d$. Denote  $\gsf = \det g$.
The principal part of the elliptic operator $\Ellop$ is
\begin{align}
  \label{eq:princip_part_normal_geo_coord}
  \Ellop_0=\gsf^{-1/2} D_d\,  \gsf^{1/2} D_d + R(x,D_{x}')
  \ \ \avec \ \
  R(x,D_{x}') = \gsf^{-1/2} \sum_{1\leq i, j\leq d-1} \! \! D_i\, \gsf^{1/2}\, g^{\prime ij} D_j,
\end{align}
with the associated principal symbol 
\begin{align}
  \label{eq:conjugated_symbol}
  \ellop(x,\xi) = \xi_d^2 + r(x,\xi') \ \ \avec \ \ r(x,\xi') =
  \sum_{1\leq i, j\leq d-1} g^{\prime ij} \xi_i \xi_j,
\end{align}
where $(g^{\prime ij}) = (g_{ij}')^{-1}$. One has $r(x,\xi') \gtrsim |\xi'|^2$. 
Note that $\Ellop_0$ is a symmetric operator. Denote by $\tilde{r}(x,\xi', \eta')$ the bilinear form in $(\xi',\eta')$
associated with the quadratic form $r(x,\xi')$.

Since the Carleman estimates we derive are `insensitive' to the
addition of any lower-order terms, it is sufficient to prove them for
$\Ellop$ replaced by $\Ellop_0$. We thus use the notation $\Ellop$ for
the operator in \eqref{eq:princip_part_normal_geo_coord} by abuse of
notation.
Observe that in the chosen local coordinates, one has
$\dnu=-\d_{x_d}$. In what follows,
denote
\begin{align*}
  P_\sigma:=P-\sigma^2  \ \ \et \ \ B = i \dnu = D_d.
\end{align*}

As is done classically,  introduce the conjugated operators
\begin{align*}
  \Psigphi=e^{\tau \varphi} \Psig e^{-\tau \varphi}
  \ \ \et \ \
  \Bphi =e^{\tau \varphi\bt} B e^{-\tau \varphi \bt }.
\end{align*}
Their respective principal symbols are
\begin{align*}
  \psigphi= \ellop\big(x,\xi+ i \tau d \varphi(x)\big) - \sigma^2 
  \ \ \et \ \
  \bphi = \xi_d + i \tau \d_d \varphi, 
\end{align*}
observing that $e^{\tau \varphi} D_j e^{-\tau
  \varphi}=D_j+i \tau \d_{j} \varphi$.
For a function $u$ as in the statements of Theorem~\ref{th:Carleman_boundary1} and Theorem~\ref{th:Carleman_boundary2}, set $v=e^{\tau \varphi} u$. Then, $e^{\tau \varphi}
\Psig u=\Psigphi v$ and $e^{\tau \varphi \bt} B u
\br=\Bphi v \br$.

Introduce the following formally
self-adjoint operators
\begin{align*}
  P_2=\frac12 \big(\Psigphi + (\Psigphi)^*\big)
  \ \ \et \ \
  P_1=\frac{1}{2 i}\big(\Psigphi - (\Psigphi)^*\big).
\end{align*}
Then $\Psigphi=P_2+i P_1$. One has $P_2 = \Psig - p(x, \tau d_x \varphi)$ and 
\begin{align*}
  P_1 = \tau  \big( \gsf^{-1/2} D_d\, \gsf^{1/2} (\d_{d} \varphi u)+ \d_{d} \varphi D_d \big)
  +\tau \sum_{j, k=1}^{d-1}
  \big(
  \gsf^{-1/2}  D_i (\gsf^{1/2}\, g^{\prime ij} (\d_j \varphi) )
  + (\d_i \varphi)  g^{\prime ij} D_j
  \big).
\end{align*}
Their respective principal symbols are
\begin{align}
  \label{eq:symbols_p2_p1}
  p_2(x,\xi,\tau,\sigma) = \ellop(x,\xi) - \sigma^2 - \ellop(x,\tau d\varphi(x))
  \ \ \et \ \
  p_1(x,\xi,\tau,\sigma) = \tellop (x, \xi, \tau d \varphi(x)),
\end{align}
where $\tellop(x,\xi,\eta)$ is the bilinear form in $(\xi,\eta)$
associated with the quadratic form $\ellop(x,\xi)$.

\begin{lemma}[\se]
  \label{lemma:sub-ellipticity}
  Suppose $\psi \in \Cinf(\M)$ is such that $|d\psi|\neq 0$ in
  $\ovl{W}$ and set $\varphi = \exp(\csp \psi)$ with $\csp >0$.  Suppose $C_0>0$. Then,
  there exists $\csp_0>0$ such that  the \se property of
  Definition~\ref{def:sub_ellipticity} holds for $(\Psig, \varphi)$ if
  $\tau \geq C_0 |\sigma|$ and $\csp \geq \csp_0$. Moreover, there
  exist $\kappa>0$ and $C_1>0$ such that
  \begin{align}
    \label{eq:sub-ellipticity-Psig}
    \kappa |\psigphi|^2 + \frac{\tau}{2i} \{ \ovl{\psigphi}, \psigphi\} \geq C_1 \lambda^4
    \qquad \si \ \ \tau \geq C_0 |\sigma|
    \ \ \et \ \ \csp \geq \csp_0,  \ \ \sigma \in \R. 
  \end{align}
\end{lemma}
\begin{proof}
  Set $\ppsi(x,\xi,\tau) = p(x,\xi + i \tau d\psi(x))$.
  By Lemmata 3.4, 3.5, and 3.55 (and their proofs) in \cite{JLR1},
  one has
  \begin{align*}
    \frac{1}{2i}\{ \ovl{\psigphi}, \psigphi\}(x,\xi, \tau)
    &= \frac{1}{2i}\{ \ovl{\pphi}, \pphi\}(x,\xi, \tau)
    = \tau H_p^2 \varphi (x,\xi) + \tau^3 H_p^2 \varphi (x,d \varphi)\\
    &=\tau \csp^4 \varphi^3 \big( (H_p \psi(x,\eta)\big)^2 + 4 \tau^2 p(x, d\psi(x))^2\big)
    + \frac{(\csp \varphi)^3}{2i}  \{ \ovl{\ppsi}, \ppsi\}(x,\eta, \tau)\\
    &\geq  C \csp \ttau^3 - C'(\ttau^3 + |\xi|^3),
  \end{align*}
  where $C$ and $C'$ are positive constants, $\ttau = \tau \csp
  \varphi$, and $\eta = (\csp \varphi)^{-1}\xi$.
Observe that $\Re \psigphi = p_2=0$ gives $|\xi|^2
\eqsim \ttau^2 + \sigma^2$. Thus, for $\csp$ chosen \suff large, the \se property holds.

Suppose  $C_0>0$ and define the compact set
\begin{align*}
  \mathbb S = \big\{ (x,\xi, \tau, \sigma)\in \ovl{W}\times \R^d \times [0,+\infty[ \times \R; \
          |\xi|^2 + \tau^2 + \sigma^2 =1 \  \ \et \ \
          \tau \geq C_0 |\sigma| \big\}.
\end{align*}
From the \se property, there exists $\kappa>0$ such that
$\kappa |\psigphi|^2 + \frac{\tau}{2i} \{ \ovl{\psigphi}, \psigphi\} \gtrsim 1$ in $\mathbb S$. Therefore, \eqref{eq:sub-ellipticity-Psig} follows by homogeneity of degree 4 in $(\xi, \tau, \sigma)$.
\end{proof}

\subsection{A first estimate: proof of Theorem~\ref{th:Carleman_boundary2}}
Here, we sketch the proof of Theorem~\ref{th:Carleman_boundary2}. The
computations are also of use for the (more delicate) proof of
Theorem~\ref{th:Carleman_boundary1}. First, one classically writes
\begin{align}
  \label{eq:Carleman_method}
  \Norm{\Psigphi v}{L^2}^2
  = Q(v) +  2 \Re \inp{P_2 w}{iP_1 w}_+,
  \qquad Q(v) = \Norm{P_2 v}{+}^2 + \Norm{P_1 v}{L^2}^2.
\end{align}
Integrations by parts as  in
(3.4.19)-(3.4.20) in \cite{JLR1}, adapted to the presence of the
spectral parameter in the term $-\sigma^2$ in $P_2$ give
\begin{align}
  \label{eq:commutateur_termes_bord}
  2 \Re \inp{P_2 v}{i P_1 v}_+
  =\Re \biginp{i [P_2,P_1] v}{v}_+ + 2 \tau \Re \tilde{\B}(v),
\end{align}
with 
\begin{align}
  \label{eq: boundary quadratic form Carleman Lopatinskii}
  \tilde{\B}(v)
  &=  
  \inp{\tilde{r} (x,D', d_{x'}\varphi) v\br}{D_dv\br}_\d
  +  \biginp{ \tilde{r} (x,D', d_{x'} \varphi) D_d v\br}{v\br}_\d\\
  & \quad
  +\inp{\d_d\varphi D_d  v\br}{D_d v\br}_\d 
  + \biginp{  \d_d\varphi \big(p(x,\tau d\varphi)+ \sigma^2 -\Opt(r)\big)
 v\br}{v\br}_\d  \notag\\
 & \quad
 +  \inp{\Op(c_0) v\br}{ D_dv\br}_\d 
 +  \biginp{  \Op(\tilde{c}_0)  D_d+  \Op(c_1)\big)  v\br}{v\br}_\d,\notag
\end{align}
with $\Op(c_0), \Op(\tilde{c}_0) \in \D^0$ and $\Op(c_1) \in \Dsct^1$.
One has $|\tilde{\B}(v)| \lesssim \normsc{\trace(v)}{1,0}^2$ as $|\sigma| \leq C_0^{-1} \tau$.
Write
\begin{align*}
  \Norm{\Psigphi v}{L^2}^2
  &= Q(v) +  \Re \biginp{i[P_2,P_1] w}{w}_+
  + 2 \tau \Re \tilde{\B}(w)\\
  &\geq \tau^{-1} \big(\kappa  Q(v) +  \tau \Re \biginp{i[P_2,P_1] w}{w}_+\big)
  + 2 \tau \Re \tilde{\B}(w),
\end{align*}
with $\kappa$ as given by Lemma~\ref{lemma:sub-ellipticity} and $\tau \geq \kappa$.
The quantified \se property given in the
second part of Lemma~\ref{lemma:sub-ellipticity} and the G{\aa}rding inequality of Theorem~6.17 in \cite{JLR2} for interior quadratic forms of type $(2,0)$ yields
\begin{align*}
  \Normsc{v}{2}^2
  \lesssim \kappa Q(v)  +  \Re \biginp{i[P_2,P_1] w}{w}_+
  + \normsc{\trace(v)}{1,1/2}^2, 
\end{align*}
for $\sigma \in \R$, and for $\tau \geq C_0 |\sigma|$ chosen \suff large. 
This yields the estimate \eqref{eq:Carleman_boundary2}, thereby proving Theorem~\ref{th:Carleman_boundary2}.\hfill \qedsymbol \endproof

\subsection{Symbol factorization and microlocal estimates}
\label{sec:microlocal_estimates}
To prove Theorem~\ref{th:Carleman_boundary1}, assume that $\d_d
\varphi >0$ at $\d\M$. If needed, reduce $W$ so that $\d_d \varphi>0$
in $\ovl{W}$.  By \eqref{eq:conjugated_symbol}, write $\psigphi
(x,\xi,\tau,\sigma) = (\xi_d + i \tau \d_d \varphi)^2 + r(x,\xi' + i
\tau d\varphi(x)) - \sigma^2$.  Introduce $\beta(x,\xi',\tau,\sigma)$
given by $\Re \beta \geq 0$ and $\beta^2 = r(x,\xi' + i \tau
d\varphi(x)) - \sigma^2$.  Then,
\begin{align*}
  \psigphi  = (\xi_d - \pi^+)(\xi_d - \pi^-), \qquad \avec \ \
  \pi^\pm(x,\xi',\tau,\sigma) = - i \tau \d_d \varphi(x) \pm i \beta(x,\xi',\tau,\sigma).
\end{align*}
The roots $\pi^\pm$ are homogeneous of degree 1 in $(\xi,\tau,\sigma)$. They are continuous functions and are smooth away from crossings. 
Note that $\Im \pi^- < - \tau \d_d \varphi <0$.

\begin{lemma}
  \label{lemma:alpha=0}
  Suppose $x\in W$, $(\xi, \tau, \sigma) \neq (0,0,0)$, and $\tau \geq
  C_0 |\sigma|$ for some $C_0>0$. If $\Re \beta =0$ then $\Im \pi^+ < 0$.
\end{lemma}
\begin{proof}
  If $\Re \beta=0$ then $\beta \in i \R$ and $\beta^2 \in \R^-$. This implies
  $r(x,\xi')\leq r(x,\tau d\varphi(x)) + \sigma^2 \lesssim \tau^2$. Thus $\tau\neq 0$
  and, as a result,  $\Im \pi^+ = - i \tau \d_d \varphi(x) < 0$. 
\end{proof}

The estimate of Theorem~\ref{th:Carleman_boundary1} is first proven
microlocally. The proof scheme is based on the sign of $\Im \pi^+$.
Set $\Ddphi = e^{\tau \varphi} D_d e^{-\tau \varphi} = D_d + i \tau \d_d \varphi$. 
\begin{proposition}
  \label{prop:microlocal_Carleman}
  Suppose $(\Psigphi, \varphi)$ has the \se property in $W$, $\dnu
  \varphi<0$ in $W \cap \d\M$, and $C_0>0$.  Suppose $x^0 \in
  W$, $\xi^{0\prime} \in \R^{d-1}$, and $\tau^0 \geq 0$ and
  $\sigma^0\geq 0$ are such that $(\xi^{0\prime}, \tau^0,\sigma^0)
  \neq (0,0,0)$ and $\tau^0 \geq  C_0 \sigma^0$. There exists a conic
  open \nhd $\U$ of $(x^0, \xi^{0\prime}, \tau^0, \sigma^0)$ in $W
  \times \R^{d-1} \times \R_+\times \R_+$ in which $\tau \geq
  C_0|\sigma|/2$ such that for $\chi \in \Ssct^0$, homogeneous of degree
  $0$ in $(\xi',\tau,\sigma)$, with $\supp \chi  \subset \U$, there
  exist $C>0$ and $\tau_*>0$ \st
  \begin{align}
    \label{eq:microlocal_Carleman}
    \tau^{-1/2}\Normsc{\Opt(\chi)v}{2}
    +\normsc{\trace(\Opt(\chi)v)}{1,1/2}
    \leq C \Big( 
    \Norm{\Psigphi v}{+}
    + \normsc{\Ddphi v\br}{1/2}
    + \Normsc{v}{2,-1}
    \Big),
  \end{align}
  for $\tau \geq \tau_*$ and $v\in \ovl{\S}(\Rdp)$.
\end{proposition}
Note that because of the condition imposed on the \nhd $\U$, one imposes in fact $\tau \geq \max (\tau_*, C_0 |\sigma|)$.

\begin{proof}[\bfseries Proof of Theorem~\ref{th:Carleman_boundary1}]
  With 
  Proposition~\ref{prop:microlocal_Carleman} at hand, one can adapt
  {\em mutatis mutandis} the proof of Proposition~8.22 in \cite{JLR2}, patch the microlocal estimates together,
  and obtain the estimate \eqref{eq:Carleman_boundary1}, which yields the result of Theorem~\ref{th:Carleman_boundary1}.
\end{proof}

\begin{proof}[\bfseries Proof of Proposition~\ref{prop:microlocal_Carleman}]
  Suppose $(x^0, \xi^{0\prime}, \tau^0, \sigma^0)$ is as in the statement.
  Two  mutually exclusive cases are  considered:
  first, the case $\Im \pi^+ (x^0, \xi^{0\prime},\tau^0,\sigma^0)\geq
    0$ and, second, the case $\Im \pi^+ (x^0, \xi^{0\prime},\tau^0,\sigma^0)<0$.

\subsubsection*{\bfseries Case 1: One root in the upper complex half-plane.}
Assume $\Im \pi^+ (x^0, \xi^{0\prime},\tau^0,\sigma^0)\geq
0$. Then, $\Re \beta(x^0, \xi^{0\prime},\tau^0,\sigma^0)>0$ by
Lemma~\ref{lemma:alpha=0}.  Since $\Im \pi^- (x^0,
\xi^{0\prime},\tau^0,\sigma^0)< 0$ and because of the homogeneity of
the roots and their continuity, there exist $C>0$ and a conic \nhd
$\U$ of $(x^0, \xi^{0\prime},\tau^0,\sigma^0)$ where $\tau \geq C_0
|\sigma|$ such that $\pi^+ (x, \xi',\tau,\sigma) \neq \pi^-(x,
\xi',\tau,\sigma)$ and
\begin{align}
  \label{eq:case1}
  \Im \pi^+ (x, \xi',\tau,\sigma)\geq -C \lambda, 
  \qquad 
  \Im \pi^- (x, \xi',\tau,\sigma)\leq -C \lambda,
  \ \ \et \ \
  \Re \beta(x, \xi',\tau,\sigma) \geq C\lambda, 
\end{align}
if $(x, \xi',\tau,\sigma) \in \U$.
Thus, $\pi^\pm \in \Ssct^1$ microlocally in $\U$.

Suppose $\chi \in \Ssct^0$ is as in the statement. Consider also $\tchi
\in \Ssct^0$, homogeneous of degree 0 and such that $\supp \tchi
\subset \U$ and $\tchi \equiv 1$ on a \nhd of $\supp \chi$. Then,
$\tpi^\pm = \tchi \pi^\pm$ and $\tbeta =\tchi\beta$ are in  $\Ssct^1$. Set $P^\pm= D_d -\Opt(\tpi^\pm) = \Ddphi \mp i \tbeta$.

Set $w = \opt(\chi) v$. By the ellipticity of $\tbeta$ in a \nhd of
$\supp \chi$, by means of a microlocal parametrix, one obtains the
following estimation (see Proposition~2.37 in \cite{JLR1}), for any
$N\in \N$,
\begin{align*}
  \normsc{w\br}{3/2} \lesssim \normsc{\opt(\tbeta)w\br}{1/2} + \normsc{v\br}{-N},
\end{align*}
for $\tau$ chosen \suff large. 
One deduces 
\begin{align*}
  \normsc{w\br}{3/2}
  \lesssim \normsc{\Ddphi w\br}{1/2}
  + \normsc{P^+ w\br}{1/2} + \normsc{v\br}{-N}.
\end{align*}
Since $\normsc{D_d w\br}{1/2} \lesssim \normsc{\Ddphi w\br}{1/2} + \normsc{w\br}{3/2}$ one obtains
\begin{align}
  \label{eq:Lopatinskii_dominate_traces}
  \normsc{\trace(w)}{1,1/2}
  &\lesssim  
    \normsc{\Ddphi  w\br}{1/2}
    + \normsc{P^+ w\br}{1/2}
    +\normsc{v\br}{-N}.
\end{align}

By \eqref{eq:case1} and Lemma~6.20 in \cite{JLR2}, for $\tau$
chosen \suff large, one has 
\begin{align}
  \label{eq: elliptic estimate case 1}
    \Normsc{\Opt(\chi) u}{1}
    + \normsc{\Opt(\chi) u\br}{1/2}
    \lesssim \Norm{P^- \Opt(\chi) u}{+}
    + \Normsc{u}{0,-N}.
\end{align}
for $u \in \ovl{\S}(\Rdp)$. 
Apply estimate~\eqref{eq: elliptic estimate case 1} to
$u = P^+ v$. It yields 
\begin{align*}
    \normsc{\Opt(\chi) P^+ v\br}{1/2}
    &\lesssim \Norm{P^- \Opt(\chi) P^+ v}{+}
    +\Normsc{v}{1,-N}
    \lesssim\Norm{\Psigphi v}{+}
    + \Normsc{v}{1},
 \end{align*}
using that $P^- \Opt(\chi) P^+ 
=\Opt(\chi)  \Psigphi  \mod \Psiscl^{1,0}$. 
Observe also that  
\begin{align*}
  \normsc{ P^+ w\br}{1/2} 
  &\lesssim \normsc{\Opt(\chi) P^+ v\br}{1/2} 
  + \normsc{ v\br}{1/2} 
  \lesssim \normsc{\Opt(\chi) P^+ v\br}{1/2} 
  + \Normsc{ v}{1},
\end{align*}
using the trace inequality of Proposition~\ref{prop:trace_inequality}. We thus obtain 
\begin{align*}
    \normsc{ P^+ w\br}{1/2} 
    &\lesssim \Norm{\Psigphi v}{+}
    +\Normsc{v}{1}.
 \end{align*}
Together with \eqref{eq:Lopatinskii_dominate_traces}, using the trace inequality of Proposition~\ref{prop:trace_inequality}, this yields the
estimate
\begin{align}
  \label{eq: Lopatinskii case 1 - trace estimation}
    \normsc{\trace(w)}{1,1/2} \lesssim  
    \normsc{\Ddphi w\br}{1/2}
    +  \Norm{\Psigphi v}{+}
    + \Normsc{v}{1},
  \end{align}
By Theorem~\ref{th:Carleman_boundary2}, as $[\Psigphi,\Opt(\chi)] \in
\Psisc^{1,0}$, one has  (see Propositions \ref{prop:calculus} and \ref{prop:operator_bounds}),
\begin{align}
  \label{eq: Lopatinskii Carleman square final estimate bis}
  \tau^{-1/2} \Normsc{w}{2}
  \lesssim \Norm{\Psigphi v}{+} 
  + \Normsc{v}{1} 
  +\normsc{\trace(w)}{1,1/2}.
\end{align}
Finally, by \eqref{eq: Lopatinskii case 1 - trace
  estimation} and \eqref{eq: Lopatinskii Carleman square final estimate bis}, 
one obtains 
\begin{align*}
    \tau^{-1/2}\Normsc{w}{2}
    +\normsc{\trace(w)}{1,1/2}
  &\lesssim \Norm{\Psigphi v}{+} + \normsc{\trace(w)}{1,1/2} + \Normsc{ v}{1} \\
  &\lesssim
    \Norm{\Psigphi v}{+}
    + \normsc{\Ddphi w\br}{1/2}
    + \Normsc{v}{1}.
  \end{align*}
By a commutator argument and the trace inequality of Proposition~\ref{prop:trace_inequality}, one obtains
\begin{align*}
  \normsc{\Ddphi w\br}{1/2} 
  &\lesssim  \normsc{\Ddphi v\br}{1/2}  
  +  \normsc{v\br}{1/2}
  \lesssim  \normsc{\Ddphi v\br}{1/2}  
  +  \Normsc{v}{1}.
\end{align*}
We then obtain~\eqref{eq:microlocal_Carleman}, which concludes the proof for Case~1. 
\hfill \qedsymbol

\subsubsection*{\bfseries Case 2: Both roots in the lower complex half-plane.}
Assume
$\Im \pi^+ (x^0, \xi^{0\prime},\tau^0,\sigma^0)< 0$. Because of the homogeneity of the roots and their continuity
, there exists $C>0$ and a conic \nhd $\U$ of $(x^0, \xi^{0\prime},\tau^0,\sigma^0)$ where $\tau \geq C_0 |\sigma|$  such that
\begin{align*}
  \Im \pi^+ (x, \xi',\tau,\sigma)\leq -C \lambda 
  \ \ \et \ \ 
  \Im \pi^- (x, \xi',\tau,\sigma)\leq -C \lambda
  \qquad \si \ \  (x, \xi',\tau,\sigma) \in \U.
\end{align*}
Suppose $\chi \in \Ssct^0$ is as in the statement and set $w = \opt(\chi) v$.
For $\eps \in (0,1)$ to be set below, as
$Q(w) \geq 0$, from \eqref{eq:Carleman_method} write
\begin{align}
  \label{eq: methode Carleman Lopatinskii - racines en bas}
  \Norm{\Psigphi  w}{+}^2 \geq \eps Q (w) + 2 \Re \inp{P_2 w}{i P_1 w}_+.
\end{align}
Introduce the compact set
$\mathbb S = \{ (x,\xi,\tau,\sigma); (x,\xi',\tau,\sigma) \in \ovl{\U}, \; 
  \xi_d \in \R, \ |\xi|^2 + \tau^2 + \sigma^2 = 1 \}$.
One has $|\psigphi(x, \xi,\tau,\sigma) | \gtrsim 1$ on $\mathbb S$.  By homogeneity, one finds
$|\psigphi(x, \xi,\tau,\sigma)  \gtrsim \lsc^2$ for $(x,\xi',\tau,\sigma) \in \U$ and $\xi_d
  \in \R$. 
Recalling $p_1$ and $p_2$ from \eqref{eq:symbols_p2_p1}, one has 
\begin{align*}
  q(x, \xi,\tau) = p_2^2(x, \xi,\tau) + p_1^2(x,\xi,\tau) \gtrsim
  \lsc^4, 
  \qquad (x,\xi',\tau) \in \U \ \ \et \ \  \xi_d\in \R. 
\end{align*}
As the interior quadratic form $Q$ has $q$ for its principal symbol in
the sense of Definition~6.13 in \cite{JLR2}, the G{\aa}rding inequality of
Theorem~6.17 in \cite{JLR2} yields, for $\tau>0$ chosen \suff large, 
\begin{align}
  \label{eq: positivity interior quadratic form Q(w) Lopatinskii}
  \Re Q(w)
    &\geq
    C\Normsc{w}{2}^2
    -C' \normsc{\trace(w)}{1,1/2}^2
    - C_N \Normsc{v}{2,-N}^2,
\end{align}
Since $[P_2,P_1]  \in \tau \Dsc^2$ one has 
\begin{align}
  \label{eq:estimation_commutator}
  |\Re \biginp{i [P_2,P_1] w}{w}_+ | 
  \lesssim \tau \Normsc{w}{2,-1}^2 
  \lesssim  \tau^{-1} \Normsc{w}{2}^2 ,
\end{align}
(see Remark~\ref{remark: tau parameter suff large ch6}).
\begin{lemma}
  \label{lemma: Lopatinskii case 2 - positivity boundary quadratic
    form}
  There exist $C, C_N>0$, and $\tau_\ast>0$ such that 
  \begin{align*}
  \tau \Re \tilde{\B}(\opt(\chi) v) &\geq C \normsc{ \trace(\opt(\chi) v)}{1,1/2}^2 
  - C_N\normsc{ \trace (v)}{1,-N}^2,
\end{align*}
for $\tau \geq \tau_\ast + C_0 |\sigma|$ and $v \in \ovl{\S}(\Rdp)$. 
\end{lemma}
A proof is given below.  Note that $C_0>0$ is the constant introduced in the statement of Theorem~\ref{th:Carleman_boundary1} and Proposition~\ref{prop:microlocal_Carleman}.
With this lemma, by \eqref{eq:commutateur_termes_bord} and \eqref{eq:estimation_commutator}, one obtains
\begin{align}
  \label{commutateur et termes de bord - ter}
 2 \Re \inp{P_2 w}{i P_1 w}_+
 &\geq C \normsc{ \trace(w)}{1,1/2}^2 
 - C' \tau^{-1} \Normsc{w}{2}^2
  - C_N\normsc{ \trace (v)}{1,-N}^2\\
  &\geq C \normsc{ \trace(w)}{1,1/2}^2 
    - C' \tau^{-1} \Normsc{w}{2}^2
    - C_N' \Normsc{v}{2,-N}^2,\notag
\end{align}
using the trace inequality of Proposition~\ref{prop:trace_inequality}.
Then, \eqref{eq: methode Carleman Lopatinskii - racines en bas},
\eqref{eq: positivity interior quadratic form Q(w) Lopatinskii} and 
\eqref{commutateur et termes de bord - ter} give
\begin{align*}
  C \eps \Normsc{w}{2}^2 
  + (C - C'\eps) \normsc{\trace(w)}{1,1/2}^2
  \leq C' \big( \Norm{\Psigphi  w}{+}^2 
  +  \tau \Normsc{w}{2,-1}^2\big) 
  + C_N \Normsc{v}{2,-N}^2.
\end{align*}
Choosing $\tau>0$ \suff large and $\eps>0$ \suff small  yields, for any $N \in \N$,  
\begin{align*}
  \Normsc{w}{2}^2 + \normsc{ \trace(w)}{1,1/2}^2
  \lesssim \Norm{\Psigphi w}{+}^2 
  + C_N \Normsc{v}{2,-N}^2 \lesssim \Norm{\Psigphi v}{+}^2
  + \Normsc{v}{1}^2
  + C_N \Normsc{v}{2,-N}^2,
\end{align*}
since $[\Psigphi , \opt(\chi)] \in \Psisc^{1,0}$. This 
yields
\begin{align*}
  \Normsc{w}{2} + \normsc{ \trace(w)}{1,1/2}
  \lesssim \Norm{\Psigphi  v}{+} 
  +\Normsc{v}{2,-1},
\end{align*}
which concludes the proof for Case~2 of the proof of Proposition~\ref{prop:microlocal_Carleman}.
\end{proof}

\begin{proof}[\bfseries Proof of Lemma~\ref{lemma: Lopatinskii case 2 -
      positivity boundary quadratic form}]
  Suppose $\tchi \in \Ssct^0$ has the same
  properties as $\chi$, that is, $\tchi$ is homogeneous
  of degree $0$ and $\supp \tchi \subset \U$, and suppose in addition
  that $\tchi =1$ on a \nhd of $\supp \chi$.
  Set $W = (\Lsct^1 w\br, D_dw\br)$  and $\tW = \Opt(\tchi) W$. By \eqref{eq: boundary quadratic
    form Carleman Lopatinskii},  write
\begin{align}
  \label{eq: positivity boundary quadratic form Lopatinskii}
  \Re \tilde{\B}(w) &\geq  \Re  \inp{\Opt(a) W}{W}_{\d} 
  - C \normsc{ \trace (w)\br}{1,-1/2}^2\\
  &\geq  \Re  \inp{\Opt(a) \tW}{\tW}_{\d} 
  - C \normsc{ \trace (w)\br}{1,-1/2}^2 -C_N \normsc{ \trace (v)}{1,-N}^2,\notag
\end{align}
with $a = (a^{ij})_{1\leq i ,j\leq 2} \in \Ssct^0$ given by
\begin{align*}
  &a^{11}=  \lsct^{-2} \d_d\varphi\br 
    \big(p(x,\tau d\varphi) + \sigma^2- r(x,\xi') \big)\br, \quad a^{22} = \d_d\varphi\br, \\
  &a^{12} = a^{21}=  \lsct^{-1} \tilde{r} (x,\xi', d_{x'}\varphi)\br,
\end{align*}
using that $|\sigma| \lesssim \tau$ to characterize some remainder terms.  
As  $\d_d\varphi\br>0$,  for $Z \in \R^2$ one has
  \begin{align}
    \label{eq: symbol positivity boundary quadratic form Lopatinskii}
    \biginp{a (x,\xi',\tau) Z}{Z}_{\R^2} \geq C  \Norm{Z}{\R^2}^2, \quad (x,\xi',\tau) \in \U, 
    \ \tau \in [1,+\infty), \quad |(\xi',\tau)| \geq R,
    \end{align}
for some $C$ and $R$ positive if and only if $\det a (x,\xi',\tau) >0$
in $\U$, by compactness
and homogeneity.
Compute 
\begin{align*}
  \lsct^{2} \det a (x,\xi',\tau)=  (\d_d\varphi)^2
  \big(p(x,\tau d\varphi) + \sigma^2- r(x,\xi')\big)\br - \tilde{r}^2 (x,\xi', d_{x'}\varphi)\br.
\end{align*}
In $\U$, as $\Im \pi^+ <0$ then $0\leq \Re \beta < \tau \d_d \varphi$.  By Lemma~\ref{lemma:prep_root_behaviour} below,
 and as $\Re \beta^2 =  r(x,\xi') -  r(x,\tau d_{x'}\varphi(x)) - \sigma^2$ and $\Im \beta^2 = 2 \tilde{r}(x,\xi',\tau d_{x'}\varphi(x))$,  this 
reads
\begin{align}
  \label{eq:case2_det>0}
   (\tau \d_d \varphi)^2 \big( r(x,\xi') - \sigma^2 -  r(x,\tau d_{x'}\varphi)\big) 
  -   (\tau \d_d \varphi)^4 + \tilde{r}^2 (x,\xi',\tau d_{x'}\varphi) < 0.
\end{align}
This is precisely $\det a (x,\xi',\tau) >0$. Note that \eqref{eq:case2_det>0} can be written, using that $\sigma \lesssim \tau$, 
\begin{align*}
  (\d_d \varphi)^2 r(x,\xi') 
  + \tilde{r}^2 (x,\xi', d_{x'}\varphi)
  <  (\d_d \varphi)^2 \sigma^2 + \tau^2 ( \d_d \varphi)^4 + \tau^2 (\d_d \varphi)^2 r(x, d_{x'}\varphi)
  \lesssim \tau^2,
\end{align*}
implying that $|\xi|\lesssim \tau$ and thus  $\lsct \lesssim \tau$ in $\U$. 

By \eqref{eq: positivity boundary quadratic form Lopatinskii}, \eqref{eq: symbol positivity boundary quadratic form Lopatinskii}
and Theorem~2.32 in \cite{JLR1}, for any
$N \in \N$, we have 
\begin{align*}
  \Re \tilde{\B}(w) &\geq C \norm{ \tW}{\d}^2 
  -C' \normsc{ \trace (w)}{1,-1/2}^2- C_N \normsc{ \trace (v)}{1,-N}^2.
\end{align*}
Observe that $\normsc{\tW}{s+s'} \lesssim \tau^s \normsc{ \tW}{s'}  + \normsc{ \trace (v)}{1,-N}$, for any $N\in \N$, 
by the G{\aa}rding inequality of Theorem~2.29 in \cite{JLR1}, since $\lsct \lesssim \tau$ in $\U$ as seen above. Thus, 
\begin{align*}
  \tau \Re \tilde{\B}(w)
  &\geq C \normsc{ \tW}{1/2}^2 
  -C' \tau \normsc{\trace (w)}{1,-1/2}^2- C_N \normsc{ \trace (v)}{1,-N}^2\\
  &\geq C \normsc{\trace(w)}{1,1/2}^2
  -C'\tau^{-1} \normsc{\trace(w)}{1,1/2}^2- C_N \normsc{ \trace (v)}{1,-N}^2.
\end{align*}
The result thus follows for $\tau$ chosen \suff large. 
\end{proof}

\begin{lemma}
\label{lemma:prep_root_behaviour}
Let $z\in\C$ such that $m=z^2$. For $r\in\R$ such
that $x_0\neq 0$, we have
\begin{equation*}
|\Re z|\lesseqqgtr |r|
\iff
4r^2 \Re m-4r^4+(\Im m)^2
\lesseqqgtr0.
\end{equation*}
\end{lemma}
See \cite[Lemma~8.11]{JLR2}. The notation $\lesseqqgtr$ is meant to represent $<$, $=$, or $>$
in each term of the equivalence.

\section{Spectral properties and resolvent estimate}
\label{sec:spectral_resolvent_estimate}
In this section, we prove
Proposition~\ref{prop:resolvent_estimate}. First, we show the spectrum
of $\SgGen$ defined in \eqref{wp-4} is localized in the half complex plane $\{ \Re z <
0\}$. Second, we prove the resolvent estimate~\eqref{eq:resolvent_est}
that holds on the imaginary axis. 

\subsection{Spectral properties}

Denote by $\Sigma_p$ and $\Sigma$ the point spectrum
and the spectrum  respectively.
\begin{proposition}\label{proposition: spectrum}
  One has $\Sigma(\SgGen) = \Sigma_p(\SgGen)$.
  If $\alpha>0$ on some open set of $\Gamma_1$, 
then $\Sigma(\SgGen)\subset \{ z\in \C; \ \Re z < 0\}$.
\end{proposition}
\begin{proof}
  Prove first that $\SgGen$ has a compact inverse. Suppose $F= (f_1,
  f_2, f_3) \in \H$ and seek $U =(u,v, w)\in D(\SgGen)$, with $w=
  v_{|\Gamma_1}$, such that $\SgGen U = F$. This reads
  \begin{align*}
    v= f_1 \in H^1_{\Gamma_0}(\M),
    \qquad \Delta_g u = f_2 \in L^2(\M),
    \qquad \et\ \ \dnu u_{|\Gamma_1} =  - \alpha v_{|\Gamma_1} - f_3
    \in L^2(\Gamma_1), 
  \end{align*}
  with $u_{|\Gamma_0}=0$ in addition.  The elliptic theory yields a
  unique solution $u \in H^{3/2}(\M)\cap H^1_{\Gamma_0}$ (use for
  instance interpolation as in \cite[Theorem~9.6]{LM:68} in the result
  of Theorem 18.63 in \cite{JLR2}). One thus finds $U \in
  (H^{3/2}(\M)\cap H^1_{\Gamma_0}(\M)) \times H^1_{\Gamma_0}(\M)
  \times H^{1/2}(\Gamma_1)$ such that $\SgGen U = F$. By
  \cite[Theorem~16.1]{LM:68}, the injection of this space into $\H$ is
  compact, yielding the compactness of $\SgGen^{-1}$.

  The spectrum of $\SgGen^{-1}$ is made of a sequence of eigenvalues
  that converges to $0$; see \cite{Brezis:11} for instance.  For
  $\lambda\neq 0$, $\lambda \in \Sigma(\SgGen^{-1})$ if and only if
  $\lambda^{-1} \in \Sigma(\SgGen)$. Hence, $\SgGen$ has a pure point
  spectrum.  Thus, to prove that $z \notin \Sigma(\SgGen)$ it suffices
  to prove that $z\id_{\H}-\SgGen$ is injective.
  
   Suppose $U=(u,v,w)\in D(\SgGen)$ and $\Re z\geq 0$.  Use that $w=
  v_{|\Gamma_1}$ and $\alpha \geq 0$ to compute
  \begin{align*}
    \biginp{(z\id_{\H}-\SgGen)U}{U}_{\H}
    &= z\Norm{U}{\H}^2 - \inp{v}{u}_{H_{\Gamma_0}^1(\M)} + \inp{P u}{v}_{L^2(\M)}
    +\inp{\dnu u_{|\Gamma_1} + \alpha v_{|\Gamma_1}}{w}_{L^2(\Gamma_1)}\\
    &=z\Norm{U}{\H}^{2}+ 2i \Im \inp{\nabla_g u}{\nabla_g v}_{L^2V(\M)}
    + \Norm{\alpha^{1/2} v_{|\Gamma_1}}{L^2(\Gamma_1)}^2.
  \end{align*}
  Then,
  $\Re \biginp{(z\id_{\H}-\SgGen)U}{U}_{\H}
  =(\Re z)\Norm{U}{\H}^2
  +  \Norm{\alpha^{1/2} v_{|\Gamma_1}}{L^2(\Gamma_1)}^2$.
  As $\Re z\geq 0$, one obtains
  \begin{align}
    \label{eq:resolvent_computation}
    |\inp{(z\id_{\H}-\SgGen)U}{U}_{\H}|
    &\geq(\Re z)\Norm{U}{\H}^{2}+\Norm{\alpha^{1/2} v_{|\Gamma_1}}{L^2(\Gamma_1)}^2.
  \end{align}

  \medskip\noindent
  {\bfseries Case $\bld{\Re z>0}$.}
  By \eqref{eq:resolvent_computation} and the Cauchy-Schwarz inequality, 
  one finds
  \begin{align*}
    \Norm{(z \id_{\H} -\SgGen) U}{\H} \geq (\Re z)\, \Norm{U}{\H},\qquad U \in D(\SgGen),
  \end{align*}
  yielding the injectivity of $z \id_{\H} -\SgGen$. Thus, $z \notin \Sigma(\SgGen)$.

  \medskip\noindent {\bfseries Case $\bld{\Re z=0}$.}  Suppose $z=i
  \sigma$, with $\sigma \in \R$, and suppose $U =
  (u,v,v_{|\Gamma_1})\in D(\SgGen)$ and $U \in \ker(i \sigma \id_{\H}
  -\SgGen)$, that is,
  \begin{align*}
    i \sigma u - v =0, \qquad i \sigma v + P u=0, \quad \et \quad
    \big(\dnu u+ (i \sigma +\alpha ) v\big)_{|\Gamma_1}=0.
  \end{align*}
  Recall that  $\Bsig = \dnu - \sigma^2 + i \alpha \sigma$.
  By \eqref{eq:resolvent_computation} one finds that  $u$ satisfies 
  \begin{align*}
    (P-\sigma^2) u =0 \ \ \dans \ \ \M,
    \quad \Bsig u=0\ \ \dans\ \ \Gamma_1, \ \ \et \ \ u=0 \ \ \dans \ \ \Gamma_1\cap \{ \alpha>0\}.
  \end{align*}
  Then, the unique continuation result of
  Theorem~\ref{th:interpolation_ineq} gives $u=0$. One concludes that
  $U=0$.
\end{proof}

\subsection{Proof of the resolvent estimate}
Suppose $U = (u,v,w) \in D(\SgGen)$,  $F=(f_0, f_1, f_2) \in \H$, and $\sigma \in \R$ are such that
$(i \sigma \id_{\H} - \SgGen) U = F$. This reads
\begin{align*}
  i \sigma u - v = f_0, \quad i \sigma v + \Ellop u = f_1,\ \ \et \ \ i \sigma w + (\dnu u + \alpha v)_{|\Gamma_1} = f_2.
\end{align*}
Since $U \in D(\SgGen)$, one has $w = v_{|\Gamma_1}$, and thus
\begin{align*}
  v = i \sigma u - f_0,
  \quad (\Ellop - \sigma^2) u = f,
  \ \ \et \ \
  (\dnu u + i\sigma \alpha u - \sigma^2 u)_{|\Gamma_1} = h,
\end{align*}
where $f= f_1 + i \sigma f_0$ and $h = f_2 + (\alpha + i \sigma)
{f_0}_{|\Gamma_1}$. Computing the inner-product of  the second part with $u$
gives, after an integration by parts and using that $u_{|\Gamma_0}=0$ and $\alpha \geq 0$,
\begin{align*}
  \inp{f}{u}_{L^2(\M)} &= \inp{(\Ellop - \sigma^2) u}{u}_{L^2(\M)}
  = - \sigma^2 \Norm{u}{L^2(\M)}^2 + \Norm{\nablag u}{L^2 V(\M)}^2 - \inp{\dnu u}{u}_{L^2(\Gamma_1)}\\
  &= - \sigma^2 \Norm{u}{L^2(\M)}^2 + \Norm{\nablag u}{L^2 V(\M)}^2
  - \inp{h}{u}_{L^2(\Gamma_1)} - \sigma^2 \norm{u}{L^2(\Gamma_1)}^2 + i \sigma \norm{\alpha^{1/2} u}{L^2(\Gamma_1)}^2.
\end{align*}
Computing the imaginary part gives
$\sigma \norm{\alpha^{1/2} u}{L^2(\Gamma_1)}^2=  \Im \inp{f}{u}_{L^2(\M)} +\Im \inp{h}{u}_{L^2(\Gamma_1)}$,
yielding by the Cauchy-Schwarz inequality and a trace inequality
\begin{align*}
  |\sigma| \norm{\alpha^{1/2} u}{L^2(\Gamma_1)}^2
  &\leq \Norm{f}{L^2(\M)}\,  \Norm{u}{L^2(\M)}
  + \norm{h}{L^2(\Gamma_1)}\,  \norm{u}{L^2(\Gamma_1)}
  \lesssim \big( \Norm{f}{L^2(\M)} + \norm{h}{L^2(\Gamma_1)} \big) \Norm{u}{H^1(\M)}\\
  & \lesssim \big( \Norm{f}{L^2(\M)} + \norm{h}{L^2(\Gamma_1)} \big) \Norm{u}{H^1(\M)}.
\end{align*}
Suppose $|\sigma| \geq |\sigma_0|>0$. One thus obtains
\begin{align*}
  \norm{\alpha^{1/2} u}{L^2(\Gamma_1)}^2
  & \lesssim |\sigma|^{-1} \big( |\sigma|\, \Norm{f_0}{H^1(\M)} + \Norm{f_1}{L^2(\M)}
  + \norm{f_2}{L^2(\Gamma_1)} \big) \Norm{u}{H^1(\M)}
  \lesssim \Norm{F}{\H} \Norm{u}{H^1(\M)}.
\end{align*}
Denote by $\Gamma_\alpha$ an open subset of $\Gamma_1$ where $\alpha \geq \delta>0$ for some $\delta>0$. Then,
\begin{align*}
  \norm{u}{L^2(\Gamma_\alpha)}^2
  \lesssim \Norm{F}{\H}\,  \Norm{u}{H^1(\M)}.
\end{align*}
A key result is the following lemma.
\begin{lemma}
  \label{lemma:resolvent-est}
  There exists $C>0$ such that
  $\Norm{u}{H^1(\M)} \leq Ce^{C \sigma^2}\big(\Norm{F}{\H} + \norm{u}{L^2(\Gamma_\alpha)}\big)$.
\end{lemma}
The proof of Lemma~\ref{lemma:resolvent-est} is based on the Carleman
estimate proved in Section~\ref{sec:Carleman_estimate}. A proof is given
in Section~\ref{sec:interpolation} below.

One deduces $\Norm{u}{H^1(\M)} \lesssim e^{C_1 \sigma^2}\big(\Norm{F}{\H} + \Norm{F}{\H}^{1/2} \Norm{u}{H^1(\M)}^{1/2}\big)$ yielding 
$\Norm{u}{H^1(\M)} \lesssim  e^{C_2 \sigma^2}\Norm{F}{\H}$,
by the Young inequality. Using $v = i \sigma u - f_0$ and $w =
v_{|\Gamma_1} = i \sigma u_{|\Gamma_1} - {f_0}_{|\Gamma_1}$, one
finally obtains
$\Norm{u}{H^1(\M)} + \Norm{v}{L^2(\M)} + \Norm{w}{L^2(\Gamma_1)} \lesssim  e^{C_3 \sigma^2}\Norm{F}{\H}$ 
using a trace inequality. This is the sought resolvent estimate~\eqref{eq:resolvent_est}.

\subsection{Interpolation inequality and proof of Lemma~\ref{lemma:resolvent-est}}
\label{sec:interpolation}
Recall the definition of $\Bsig$ in \eqref{eq:Bsig}.
From the Carleman estimate of
Proposition~\ref{prop:Carleman_boundary-Bsigma-1}, following the proof
of Theorem~6.21 in \cite{JLR1} and Theorem~9.4 in \cite{JLR2}, one
obtains the following interpolation inequality. 
\begin{proposition}[interpolation initiated in the interior]
  \label{prop:interpolation_ineq_1}
  Suppose $V$ and $\omega$ are open subsets of $\M$. There exist $C>0$
  and $\delta \in (0,1)$ such that
  \begin{align*}
    \Norm{y}{H^1(V)}
    \leq Ce^{C \sigma^2}
    \Norm{y}{H^1(\M)}^{1-\delta} \big(
    \Norm{(\Ellop - \sigma^2) y}{L^2(\M)}
    + \norm{\Bsig y_{|\Gamma_1\cap V}}{L^2(\Gamma_1\cap V)}
    + \norm{y}{L^2(\omega)}
     \big)^\delta,
  \end{align*}
 for $y \in H^2(\M) \cap H^1_{\Gamma_0}(\M)$.
\end{proposition}
From the Carleman estimate of
Proposition~\ref{prop:Carleman_boundary-Bsigma-2}, following the proof
of Lemma~9.5 in \cite{JLR2}, one obtains the following interpolation
inequality.
\begin{proposition}[interpolation initiated at the boundary]
  \label{prop:interpolation_ineq_2}
  Suppose $U$ is an open subset of $\M$ such that $U \cap \Gamma_1 \neq
  \emptyset$. There exist $C>0$ and $\delta \in (0,1)$ such that
  \begin{align*}
    \Norm{y}{H^1(U)}
    \leq Ce^{C \sigma^2}
    \Norm{y}{H^1(\M)}^{1-\delta} \big(
    \Norm{(\Ellop - \sigma^2) y}{L^2(\M)}
    + \norm{\Bsig y_{|\Gamma_1\cap U}}{L^2(\Gamma_1\cap U)}
    + \norm{y_{|\Gamma_1\cap U}}{L^2(\Gamma_1 \cap U)}
     \big)^\delta,
  \end{align*}
 for $y \in H^2(\M) \cap H^1_{\Gamma_0}(\M)$.
\end{proposition}
Together, Propositions~\ref{prop:interpolation_ineq_1} and
\ref{prop:interpolation_ineq_2} imply the following result.
\begin{theorem}
  \label{th:interpolation_ineq}
  Suppose $V$ is an open subset of $\M$ and $\Gamma$ is an open subset
  of $\Gamma_1$. There exist $C>0$ and $\delta \in (0,1)$ such that
  \begin{align*}
    \Norm{y}{H^1(V)}
    \leq Ce^{C \sigma^2}
    \Norm{y}{H^1(\M)}^{1-\delta} \big(
    \Norm{(\Ellop - \sigma^2) y}{L^2(\M)}
    + \norm{\Bsig y_{|\Gamma_1}}{L^2(\Gamma_1)}
    + \norm{y_{|\Gamma}}{L^2(\Gamma)}
     \big)^\delta,
  \end{align*}
 for $y \in H^2(\M) \cap H^1_{\Gamma_0}(\M)$.
\end{theorem}
Apply this interpolation inequality to $u$ with $V = \M$ and $\Gamma
= \Gamma_{\alpha}$. It gives
\begin{align*}
    \Norm{u}{H^1(\M)}
    &\lesssim e^{C \sigma^2}
    \big(
    \Norm{(\Ellop - \sigma^2) u}{L^2(\M)}
    + \norm{h}{L^2(\Gamma_1)}
    + \norm{u}{L^2(\Gamma_\alpha)}
    \big)\\
    &\lesssim e^{C \sigma^2}
    \big(
    \Norm{f}{L^2(\M)}
    + \norm{h}{L^2(\Gamma_1)}
    + \norm{u}{L^2(\Gamma_\alpha)}
    \big),
\end{align*}
which yields the result of Lemma~\ref{lemma:resolvent-est}.

\section{Interior damping}\label{sec:interior}

In system~\eqref{wp-1} damping occurs through the boundary condition. A second model is the following one in which damping acts in the interior of the domain:
\begin{align}
  \begin{cases}
     \label{eq:interior_damping}
    \d_{t}^{2}u + \Ellop u + \alpha \d_t u =0
    \ \  \dans \ \ \M \times (0,+\infty), \\
    \dnu  u+\d_{t}^{2}u=0
    \ \ \sur \ \ \Gamma_1 \times (0,+\infty)  \ \ \et \ \
    u = 0 \ \ \sur \ \Gamma_0 \times (0,+\infty),\\
    u\ic=u_0, \quad \d_t u\ic=u_1
    \ \ \dans  \ \  \M,
  \end{cases}
\end{align}
where $\alpha \in \Con(\M)$ such that $\alpha\geq 0$. Such interior
damping has also received a lot of attention.

System~\eqref{eq:interior_damping} can be recast in the form of a semigroup equation
\begin{align}
  \label{eq:interior_damping2}
    \frac{d W}{d t}(t) =  \SgGen_{int} W(t), \qquad 
    W_{|t=0}= W_0= (u_{0}, v_{0}, {v_0}_{|\Gamma_1}),
\end{align}
with the following generator
\begin{align}
    \label{eq:interior_damping3}
    \SgGen_{int}=\begin{pmatrix} 0 & \id & 0 \\
    -\Ellop & -\alpha & 0 \\ -\dnu ( \cdot)_{|\Gamma_1} & 0 & 0	
    \end{pmatrix},
\end{align}
with domain, formally identical to that of the generator in the case of boundary damping,
\begin{align}
  \label{eq:interior_damping4}
  D(\SgGen_{int})  =\Big\{ (u,v,w) \in \H; \ 
		v \in H_{\Gamma_0}^{1}(\Omega), \ 
		-\Ellop u \in L^{2}(\Omega), \
		w= v_{|\Gamma_1}, \ 
		\dnu u_{|\Gamma_1} \in L^2(\Gamma_1)
	\Big\}.
\end{align}

Well-posedness for \eqref{eq:interior_damping2} can also be obtained for weak and strong solutions with the same Hilbert space $\H$ as in \eqref{eq:Hilbert_space_H}. The definition of the energy remains unchanged; see
\eqref{eq:def_energy}. The counterpart to
Theorem~\ref{thm:stabilization} is then the following result.
\begin{theorem}
  \label{thm:interior_stabilization}
  Suppose $\omega$ is a nonempty open subset of $\M$ and
  $\alpha \in \Con(\M)$ is such that $\alpha>0$ on
  $\omega$. Suppose $k \in \N^*$. Then, there exists $C>0$ such
  that
  \begin{align*}
    \energy(W)(t) = \frac12 \Norm{W(t)}{\H}^2
    \leq \frac{C}{\log (2+t)^k} \Norm{\SgGen_{int}^k U^0}{\H}^2,
    \qquad \si \  W_0 \in D(\SgGen_{int}^k).
\end{align*}
\end{theorem}
Similarly, one proves that the spectrum of $\SgGen_{int}$ lies in $\{ \Re z
< 0\}$.  The proof of Theorem~\ref{thm:interior_stabilization} is then
equivalent to the following resolvent estimate:
there exist $K>0$ and
$\sigma_0>0$ such that
\begin{align}
  \label{eq:resolvent_est_interior}
  \Norm{(i \sigma \id_{\H}-\SgGen_{int})^{-1}}{\mathscr{L}(\H, \H)}
  \leq K e^{K \sigma^2}, \quad \sigma \in \R, \ |\sigma| \geq \sigma_0 .
\end{align}
We sketch the proof. Suppose $U = (u,v,w) \in D(\SgGen_{int})$,  $F=(f_0, f_1, f_2) \in \H$, and $\sigma \in \R$ are such that
$(i \sigma \id_{\H} - \SgGen_{int}) U = F$. This reads
\begin{align*}
  i \sigma u - v = f_0, \quad (i \sigma +\alpha) v + \Ellop u = f_1,\ \ \et \ \ i \sigma w + \dnu u_{|\Gamma_1} = f_2.
\end{align*}
As $w = v_{|\Gamma_1}$, write
\begin{align*}
  v = i \sigma u - f_0,
  \quad (\Ellop - \sigma^2 + i \sigma \alpha) u = f,
  \ \ \et \ \
  (\dnu u  - \sigma^2 u)_{|\Gamma_1} = h,
\end{align*}
where $f= f_1 + (i \sigma +\alpha)f_0$ and $h = f_2 + i\sigma {f_0}_{|\Gamma_1}$.
Computing the inner-product of  the second part with $u$
gives, after an integration by parts and using that $u_{|\Gamma_0}=0$ and $\alpha \geq 0$,
\begin{align*}
  \inp{f}{u}_{L^2(\M)} + \inp{h}{u}_{L^2(\Gamma_1)}
  &= \inp{(\Ellop - \sigma^2+ i \sigma \alpha) u}{u}_{L^2(\M)}
  + \inp{(\dnu u  - \sigma^2 u)_{|\Gamma_1}}{u}_{L^2(\Gamma_1)}\\
  &= - \sigma^2 \Norm{u}{L^2(\M)}^2 + i \sigma \Norm{\alpha^{1/2} u}{L^2(\M)}^{2}
  + \Norm{\nablag u}{L^2 V(\M)}^2 - \sigma^2 \norm{u}{L^2(\Gamma_1)}^2.
\end{align*}
Computing the imaginary part gives
$ \sigma \Norm{\alpha^{1/2} u}{L^2(\M)}^2=  \Im \inp{f}{u}_{L^2(\M)} +\Im \inp{h}{u}_{L^2(\Gamma_1)}$,
yielding by the Cauchy-Schwarz inequality and a trace inequality
\begin{align*}
  |\sigma| \Norm{\alpha^{1/2} u}{L^2(\M)}^2
   \lesssim \big( \Norm{f}{L^2(\M)} + \norm{h}{L^2(\Gamma_1)} \big) \Norm{u}{H^1(\M)}.
\end{align*}
Denote by $\omega$ an open subset of $\M$ where $\alpha \geq \delta>0$ for some $\delta>0$. Then,
\begin{align*}
  \norm{u}{L^2(\omega)}^2
  \lesssim \Norm{F}{\H}\, \Norm{u}{H^1(\M)}.
\end{align*}
The Lemma \ref{lemma:resolvent-est_interior} stated below is a consequence of
Proposition~\ref{prop:interpolation_ineq_1} and implies the
resolvent estimate~\eqref{eq:resolvent_est_interior}.
\begin{lemma}
  \label{lemma:resolvent-est_interior}
  There exists $C>0$ such that
  $\Norm{u}{H^1(\M)} \leq Ce^{C \sigma^2}\big(\Norm{F}{\H} + \Norm{u}{L^2(\omega)}\big)$.
\end{lemma}

\section{Optimality}
\label{sec:optimality}
This section addresses the optimality of the resolvent estimate.  For
boundary damping, Proposition~\ref{prop:resolvent_estimate} states the
resolvent estimate
\begin{align*}
  \Norm{(i\sigma \id_\H- \SgGen)^{-1}}{\mathcal{L}(\H,\H)} \leq K e^{K\sigma^2}.
\end{align*}
In the case of interior damping, one has a similar estimate in
\eqref{eq:resolvent_est_interior}. It is shown in, by means of
explicit examples in what follows, that the power $2$ in the spectral
parameter $\sigma$ is optimal.  Both examples for interior and
boundary damping have the same features: one constructs quasi-modes
that tend to concentrate away from the damping region. For the
aforementioned quasi-modes, the damping barely acts, as if it were
vanishing. Since $\alpha=0$ means that the spectrum of the semigroup
generator lies on $i \R$, it yields a poor bound if the resolvent $(i
\sigma\id_\H - \SgGen)^{-1}$ acts on such quasi-modes for $\sigma$
well chosen.

Define the cut-off function $\chi \in \Con^\infty([0,1])$ satisfying
\begin{equation}
\label{eq:chi_optimality}
    0 \leq \chi \leq 1, \qquad \chi \equiv 1 \text{ on } [1/2,1] \qquad  \text{and} \qquad  \supp \chi \subset [1/4,1].
\end{equation}

\subsection{Interior damping}
\label{sec:optimality_interior}

For $d \geq 2$, consider the compact manifold $\M=\mathbb{T}^{d-1}\times [1/8,1]$ where $\mathbb{T}^{d-1}=(\mathbb{R} / 2\pi \mathbb{Z})^{d-1}$ and suppose $a \in \Cinf(\M)$ is such that
$a \geq 0$ and $\supp a \subset \mathbb{T}^{d-1}\times(1/4,1/2)$. Denote by $y$ the
variable that lies in $[1/8,1]$. Let $\Gamma_0=\mathbb{T}^{d-1} \times \{1/8\}$, $\Gamma_1=\mathbb{T}^{d-1} \times \{1\}$ and set $P = - \Delta_{\mathbb{T}^{d-1}} -
\d_y^2$.

\begin{lemma}
	\label{lemma:optimality_interior}
	There exists $C>0$ such that for any $\sigma_0>0$, there exist $\sigma\geq
	\sigma_0$ and $v \in \Cinf(\M)\cap H_{\Gamma_0}^{1}(\M)$ such that
	\begin{equation*}
		\Norm{(P-\sigma^2+ia\sigma)v}{L^2(\M)}^{2}
		\leq C e^{-\sigma^2}\Norm{(v,i\sigma v, i\sigma v_{|\Gamma_1})}{\H}^{2},
	\end{equation*}
	with $(\d_\nu v-\sigma^2v)_{|\Gamma_1} =0$.
\end{lemma}
\begin{proof}
	For $k \in \mathbb{N}$, define
	\begin{equation*}
		\phi_k(x)=(2 \pi)^{-(d-1) / 2} e^{i k x_1}.
	\end{equation*}
	Representing the torus by the cube $[0,2 \pi]^{d-1}$, it results that
	\begin{align*}
		\Norm{\phi_k}{L^2(\mathbb{T}^{d-1})}^2 & =(2 \pi)^{-(d-1)} \int_{[0,2 \pi]^{d-1}}|e^{i k x_1}|^2 d x  =(2 \pi)^{-(d-1)}(2 \pi)^{d-1}=1.
	\end{align*}
	Define $m(\sigma)=\sigma^2+\sigma^4$. Set $\sigma=((\sqrt{1+4k^2}-1)/2)^{1/2}$; it is the unique positive root of $m(\sigma)=k^2$.  Choose $k$ such that $\sigma \geq \sigma_0$.  With $\chi$ given in \eqref{eq:chi_optimality}, define
	\begin{equation*}
		v(x,y)=\chi(y)e^{\sigma^2y}\phi_k(x).
	\end{equation*}
	Since $\chi$ vanishes in a neighborhood of $y=1/8$, one has $v_{|\Gamma_0}=0$. A direct computation gives
	\begin{equation*}
		(P-\sigma^2)v=-(\chi''(y)+2\sigma^2\chi'(y))e^{\sigma^2y}\phi_k(x). 
	\end{equation*}
	As $\chi'(1)=0$, one has $\dnu v_{|\Gamma_1}=\d_{y}v_{|\Gamma_1}=\sigma^2v_{|\Gamma_1}$.
	Thus, $v$ is a solution to 
	\begin{equation*}
		\begin{cases}
			P v-\sigma^2v+i\sigma a v=f \quad \dans \ \ \M,\\
			\dnu v-\sigma^2v=0\quad \sur \ \ \mathbb{T}^{d-1}\times \{1\}
			\ \ \ \et\ \ \ v=0 \quad \sur \ \  \mathbb{T}^{d-1} \times \{1/8\},
		\end{cases}
	\end{equation*}
	where $f=(i\sigma a \chi(y) -\chi''(y)-2\sigma^2\chi'(y))e^{\sigma^2y}\phi_k(x)$. 
	Since $\supp a \subset \mathbb{T}^{d-1} \times(1/4,1/2)$ and $\supp \chi' \subset [1/4,1/2]$, one obtains
	\begin{equation}\label{resolvent estimate: interior general}
		\Norm{f}{L^2(\M)}^{2} 
		\lesssim \sigma^4\int_{1/4}^{1/2}e^{2\sigma^2y}\, dy
		\lesssim \sigma^2 e^{\sigma^2}.
	\end{equation}
	For the quasimode norm, keeping only the boundary component,
	\begin{equation}\label{resolvent estimate: interior general 2}
		\Norm{(v,i\sigma v, i\sigma v_{|\Gamma_1})}{\mathcal{H}}^2 \geq \norm{i\sigma v_{|\Gamma_1}}{L^2(\Gamma_1)}^2=\sigma^2e^{2\sigma^2}.
	\end{equation}
	Combining \eqref{resolvent estimate: interior general} and \eqref{resolvent estimate: interior general 2} yields the result.
\end{proof}

\begin{figure}[htb!]
  \begin{center}
  \resizebox{7cm}{!}{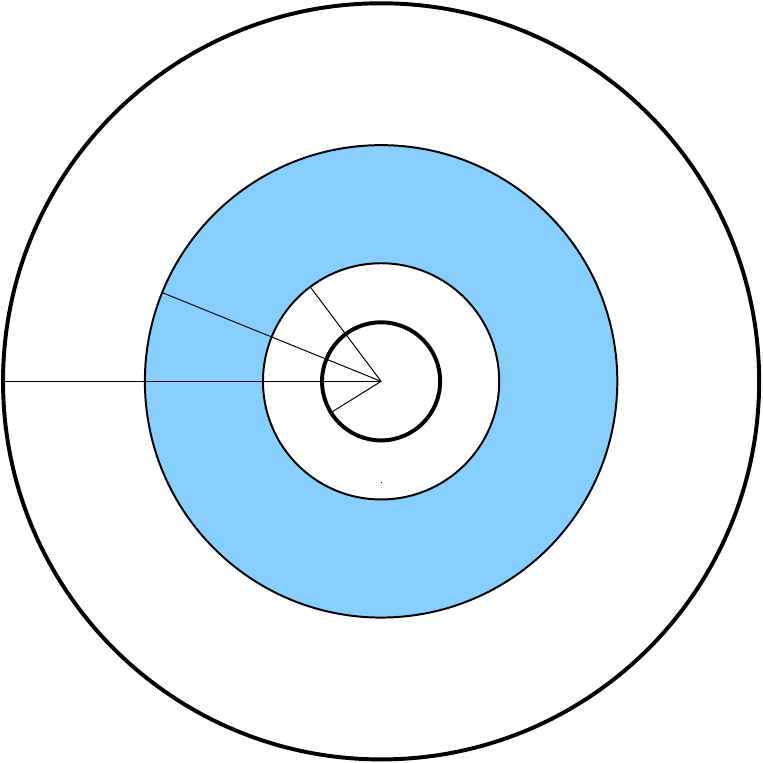}
  \caption{Sketch of an open set $\Omega\subset \R^2$ on which the result of Lemma~\ref{lemma:optimality_interior} holds. (Some proportions are not respected.)}
  \label{fig:Crown}
\end{center}
\end{figure}

For $d=2$, the example of Lemma~\ref{lemma:optimality_interior} is for $\mathcal M = \mathbb T \times [1/8,1]$,  along with the metric $g = g_{\mathbb T} + dy^2$, where $g_{\mathbb T}$ is the metric induced by the Euclidean metric on $\R^2$ on the unit circle.
This corresponds to the case of the open set $\Omega$ of $\R^2$ given by 
$\Omega = B(0,1) \setminus \ovl{B(0,1/8)}$, where $B(0,r)$ denotes the Euclidean ball centered at $0$ with radius $r$, equipped with the non-Euclidean metric given by $g = dr^2 + g_{\mathbb T}$, in polar coordinates. See Figure~\ref{fig:Crown} for a sketch.

Consider the operator
$(\SgGen_{int}, D(\SgGen_{int}))$ defined in
\eqref{eq:interior_damping3}--\eqref{eq:interior_damping4} in the present geometry. For $n\in \N$, let  $\sigma\geq n$ and $v$ be as given by  Lemma \ref{lemma:optimality_interior}. 
Set $V = (v, i \sigma v, i \sigma v_{|\Gamma_1})$ and $F= (i\sigma\id_{\H} - \SgGen_{int}) V
= (0, (P-\sigma^2+i\sigma a) v, 0)$. By \eqref{eq:resolvent_est_interior} and  Lemma \ref{lemma:optimality_interior}, one obtains
\begin{equation*}
  K e^{K\sigma^2}\geq \Norm{(i\sigma \id -\SgGen_{int})^{-1}}{\mathscr{L}(\H, \H)}
  \gtrsim  \frac{\Norm{V}{\H}}{\Norm{F}{\H}}
  \gtrsim  \frac{\Norm{V}{\H}}{\Norm{(P-\sigma^2v+i\sigma a) v}{L^2(\M)}}
  \gtrsim e^{\frac{\sigma^2}{2}}.
\end{equation*}
Since the right-hand side grows faster than $Ke^{K\sigma^\eta}$ for
any $0< \eta < 2$ and any $K > 0$, this establishes the sharpness of the
exponent~$2$.

\subsection{Boundary damping}
For the boundary damping problem, the cases $d=2$ and $d\geq 3$ are
considered separately.

\begin{figure}
  \begin{center}
  \resizebox{7cm}{!}{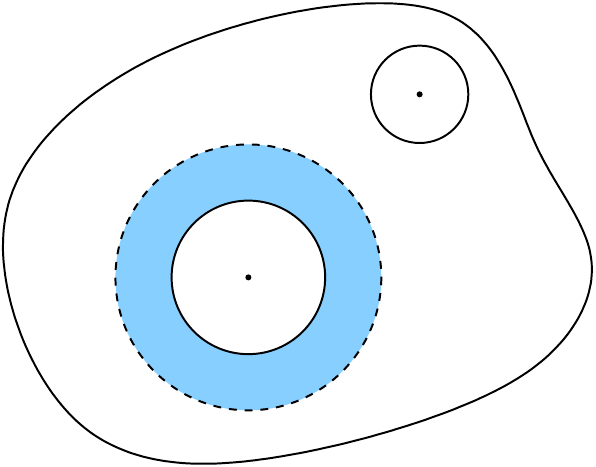}
  \caption{Geometry of the two-dimensional example with a  localized boundary damping.}
  \label{fig:example_d=2}
\end{center}
\end{figure}

\paragraph{\bfseries Case $\bld{d=2}$.}
Suppose $a \in \R^2$ and $\caO$ is a smooth open set such that
$\ovl{B(a, 2)} \subset \caO$.  Consider $b \in \caO$ and $r_b>0$ such
that $\ovl{B(b, r_b)} \subset\caO\setminus \ovl{B(a, 2)}$.  Set $\M =
\ovl{\caO} \setminus \big( B(a, 1) \cup B(b, r_b)\big)$.
One can choose $\caO$ so that $\d\M$ has three connected components:
$\Gamma_1^a = \d B(a,1)$, $\Gamma_1^b = \d B(b,r_b)$, and
$\Gamma_0$. This is illustrated in Figure~\ref{fig:example_d=2}.

Use polar coordinates  $(r, \theta)$ with $a$ as the origin. One has $r \geq 1$ in $\M$. Define on $\M$ the  smooth Riemannian metric $g = dr^2 + h(r)^2 d\theta^2$ where $h$ is a smooth function on $[1,+\infty[$ such that $h \geq 1$ and 
    \begin{align*}
      h(r) = \begin{cases}
        1 & \si \ 1\leq r\leq 7/4,\\
        r & \si \ r\geq 2.
      \end{cases}
    \end{align*}
    Set $\caU = \{ 1\leq r \leq 7/4\}$, the blue region in Figure~\ref{fig:example_d=2}. In $\caU$, $\Delta_g = \d_r^2
    + \d_\theta^2$. As $g$ is the Euclidean metric if $r\geq 2$ one
    finds $\Delta_g = \d_x^2 +\d_y^2 = \d_r^2 + r^{-1}\d_r + r^{-2}
    \d_\theta^2$ in $\M \setminus \caU$.

    Set $\Gamma_1 = \Gamma_1^a\cup \Gamma_1^b$ and define $\alpha \in \Cinf(\Gamma_1)$ by $\alpha=0$ on $\Gamma_1^a$, $\alpha=1$ on  $\Gamma_1^b$. 
\begin{lemma}
There exists $C>0$  such that for any $\sigma_0>0$ there exists $\sigma\geq \sigma_0$ and $v \in \Cinf(\M) \cap H^1_{\Gamma_0}(\M)$ such that
\begin{equation*}
  \Norm{(P-\sigma^2)v}{L^2(\M)}^{2}  \leq Ce^{-\sigma^2} \Norm{(v,i\sigma v, i\sigma v_{|\Gamma_1})}{\H}^{2},
\end{equation*}
with $(\d_\nu v  -\sigma^2 v+ i\sigma\alpha v)_{|\Gamma_1}=0$. 
\end{lemma}
\begin{proof}
Set $m(\sigma)=\sigma^2+\sigma^4$. For $k\in \mathbb{N}$, set $
\phi_k(\theta)=(2\pi)^{-1/2}e^{ik\theta}$ and
$\sigma=\big((\sqrt{1+4k^2}-1)/2\big)^{1/2}$; it is the unique positive root
of $m(\sigma)=k^2$.  Choose $k$ so that $\sigma \geq \sigma_0$.

Use $z=r-1$. Let
$\zeta \in \Cinf(\R)$ with $0\le\zeta \leq 1$,
\begin{equation*}
\zeta(z)=1 \ \text{on} \ z \leq 5/8 ,\quad \zeta(z)=0 \ \text{on} \  z\geq 11/16.
\end{equation*}
In particular, $\zeta(z)=1$ for $z\leq 1/2$ and $\zeta(z)=0$ for $z\geq 3/4$. Consequently, $\supp\zeta'\subset[1/2,3/4]$ and $\supp\zeta''\subset[1/2,3/4]$. Define
\begin{equation*}
v(\theta,z)=\begin{cases}
\zeta(z)e^{-\sigma^2z}\phi_k(\theta)
& \ \text{ on }\mathcal U,\\
0 & \ \text{ on } \mathcal{M} \backslash \mathcal{U}.
\end{cases}
\end{equation*}
Since $\zeta=0$ in a neighborhood of $z=3/4$, every derivative of $v$ vanishes there, so this extension is $\Cinf$ on $\mathcal{M}$. In particular, $v\equiv0$ near
$\Gamma_0\cup\Gamma_1^b$, hence  ${v}_{|\Gamma_0}=0$.

In $\caU$, a direct calculation gives 
\begin{align*}
	(P-\sigma^2)v&=
	((k^2-\sigma^4-\sigma^2)\zeta
  -\zeta''+2\sigma^2\zeta')e^{-\sigma^2z}\phi_k
  = (-\zeta''+2\sigma^2\zeta')e^{-\sigma^2z}\phi_k=f,
\end{align*}
using $m(\sigma)=k^2$. Hence,  $\supp f\subset [1/2,3/4] \times \mathbb{T}$.  Notice that
the coefficients of $P$ outside $\mathcal U$ have no impact on this computation.
Since $v$ vanishes in a neighborhood of $\Gamma_1^b$,
both $v$ and $\d_\nu v$ vanish there. On  $\Gamma_1^a=\{z=0\}$, one has $\zeta(0)=1$, $\zeta'(0)=0$,
and $\alpha=0$. Hence
\begin{equation*}
  \d_\nu v-\sigma^2v+ i\sigma\alpha v=0 \ \text{on} \ \Gamma_{1}.
\end{equation*}
Using the definition of $f$, one estimates
\begin{equation}\label{estimate:f}
\Norm{f}{L^2(\mathcal{M})}^2
=\int_{1/2}^{3/4}|2\sigma^2\zeta'(z)-\zeta''(z)|^2e^{-2\sigma^2z} dz
\leq C_\zeta \sigma^4\int_{1/2}^{3/4}e^{-2\sigma^2z}d z\leq C \sigma^2 e^{-\sigma^2}.
\end{equation}
Since $v_{|\Gamma_1^a}=\phi_k$, one obtains
\begin{equation}\label{estimate:quasi}
\Norm{(v,i\sigma v, i\sigma v_{|\Gamma_1})}{\mathcal{H}}^2 \geq \norm{i\sigma v_{|\Gamma_1^a}}{L^2(\Gamma_1^a)}^2=\sigma^2 \int_0^{2 \pi}\left|\phi_k(\theta)\right|^2 d \theta  =\sigma^2.
\end{equation}
The result follows by combining  \eqref{estimate:f} and \eqref{estimate:quasi}.
\end{proof}

\paragraph{\bfseries Case $\bld{d\geq 3}$.}
The construction is inspired by the appendix on spherical harmonics in \cite{Zelditch}. Let $(\mathcal{N},g)$ be a compact $(d-1)$-dimensional Riemannian manifold without
boundary, and set $\M = \mathcal{N} \times I$ with $I = [1/8,1]$, where $y$
denotes the variable lying in $I$. As before, $\Gamma_0= \mathcal{N} \times\{1 / 8\}$, $\Gamma_1=\mathcal{N} \times\{1\}$, and $P=-\Delta_{\mathcal{N}}-\d_y^2$. 

As in Section~\ref{sec:optimality_interior}, the following lemma yields the
optimality result stated in Proposition \ref{prop:optimality}. An example satisfying the
hypotheses of the lemma is given afterward.

\begin{lemma}
	\label{lem:quasimode-boundarygeneral}
	Suppose there exist $C>0$ and $C_0 >0$
	such that, for every $k_0 > 0$, there exist $k \geq k_0$ and
	$\phi_k \in H^2(\mathcal{N})$ with $-\Delta_{\mathcal{N}}\phi_k=k(k+d-2)\phi_k$ and
	$\Norm{\phi_k}{L^2(\mathcal{N})} = 1$ satisfying
	\begin{equation}\label{optimality: boundary general}
\Norm{\alpha \phi_k}{L^{2}(\mathcal{N})}^{2}
		\leq Ce^{-C_0 k}.
	\end{equation}
	Then, there exists $C'>0$  such that for any $\sigma_0>0$ there exists $\sigma\geq \sigma_0$ and $v \in \Cinf(\M) \cap H^1_{\Gamma_0}(\M)$ such that
	\begin{equation*}
		\Norm{(P-\sigma^2)v}{L^2(\M)}^{2}
		+\Norm{(\dnu v-\sigma^2v+i\alpha\sigma v)_{|\Gamma_1}}{L^{2}(\mathcal{N})}^{2}
		\leq C'  e^{-\min(C_0,1)\sigma^2} \Norm{(v,i\sigma v, i\sigma v_{|\Gamma_1})}{\H}^{2}.
	\end{equation*}
\end{lemma}
\begin{proof}
Define $m(\sigma)=\sigma^4+\sigma^2$. For $k \in \N$, set $\sigma = ((\sqrt{1+4k(k+d-2)}-1)/2)^{1/2}$; it is the unique positive root of $m(\sigma) = k(k+d-2)$. Choose $k$ such that $\sigma \geq \sigma_0$. With $\chi$ given in \eqref{eq:chi_optimality}, define
\begin{align*}
 v (x,y)= \chi(y) e^{\sigma^2 y}\phi_k(x), \qquad (x,y)\in \mathcal N\times [1/8,1].
\end{align*}
Since $\chi$ vanishes in a neighborhood of $y=1/8$, it follows that $v_{|\Gamma_0}=0$.  A direct computation gives
\begin{equation*}
	h=(P-\sigma^2) v  =-(\chi''+2 \sigma^2 \chi') e^{\sigma^2 y} \phi_k \ \text{and} \ g=(\d_y v-\sigma^2 v+i \alpha \sigma v)_{|\Gamma_1}=i \alpha \sigma e^{\sigma^2} \phi_k.
\end{equation*}

Multiplying $k(k+d-2)=\sigma^4+\sigma^2$ by $4$ and completing the squares gives $(2k+d-2)^2
= (2\sigma^2+1)^2+(d-1)(d-3)$.  Since $k\geq 0$ and $d\geq 3$, both $2k+d-2$ and $2\sigma^2+1$ are positive. Hence, $k\geq \sigma^2-(d-3)/2$.

Since $\supp \chi' \cup \supp  \chi'' \subset$ $[1 / 4,1 / 2]$ and $\Norm{\phi_k}{L^2(\mathcal{N})}=1$, one has
\begin{align*}
	\Norm{h}{L^2(M)}^2 & =\int_{1 / 4}^{1 / 2} \int_\mathcal{N}|\chi''(y)+2 \sigma^2 \chi'(y)|^2 e^{2 \sigma^2 y}\left|\phi_k(x)\right|^2 dx d y \leq C \sigma^2e^{\sigma^2},\\	\Norm{g}{L^2(\mathcal{N})}^2& \lesssim \sigma^2e^{2\sigma^2}\Norm{\alpha \phi_k}{L^2(\mathcal{N})}^2\lesssim \sigma^2e^{2\sigma^2}e^{-C_0k}\lesssim \sigma^2e^{(2-C_0)\sigma^2}.
\end{align*}

Then, it follows that
\begin{equation}\label{eq1}
	\Norm{h}{L^2(\M)}^2+\Norm{g}{L^2(\mathcal{N})}^2\lesssim\sigma^2e^{\max(1,2-C_0)\sigma^2}.
\end{equation}
On the other hand, 
\begin{equation}\label{eq2}
	\Norm{(v,i\sigma v, i\sigma v_{|\Gamma_1})}{\mathcal{H}}^{2}\geq \norm{i\sigma v_{|\Gamma_1}}{L^2(\Gamma_1)}^{2}= \sigma^2e^{2\sigma^2}.
\end{equation}
Combining \eqref{eq1} and \eqref{eq2}, the desired estimate follows.
\end{proof}

We now present an example of a manifold $(\mathcal{N},g)$ for which the hypotheses of Lemma~\ref{lem:quasimode-boundarygeneral} hold. Set $\mathcal{N}=\mathbb{S}^{d-1}$, and for $x=\left(x_1, \ldots, x_d\right) \in \R^{d}$, consider the polynomial
\begin{equation*}
	Q_k\left(x_1, \ldots, x_d\right)=\left(x_1+i x_2\right)^k
\end{equation*}
which is homogeneous of degree $k$, and define
\begin{equation*}
	\psi_k(x)={Q_k}_{|\mathbb{S}^{d-1}}, \quad \phi_k=\frac{\psi_k}{\Norm{\psi_k}{\mathbb{S}^{d-1}}}=c_k\psi_k, \quad k \in \mathbb{N}.
\end{equation*}

Computing the derivatives, one obtains $\Delta_{\mathbb{R}^d} Q_k=0 $. By the polar decomposition of the Euclidean Laplacian, $
\Delta_{\mathbb{R}^d}=\d_r^2+\frac{d-1}{r} \d_r+\frac{1}{r^2} \Delta_{\mathbb{S}^{d-1}},$ and by the homogeneity $Q_k(r \omega)=r^k Q_k(\omega)$, it follows that
\begin{equation*}
	\d_r(r^k Q_k(\omega))=k r^{k-1} Q_k(\omega), \quad \d_r^2(r^k Q_k(\omega))=k(k-1) r^{k-2} Q_k(\omega).
\end{equation*}
Therefore,
\begin{equation*}
	0=\Delta_{\mathbb{R}^d} Q_k(r \omega)=r^{k-2}\left[k(k+d-2) Q_k(\omega)+\Delta_{\mathbb{S}^{d-1}} Q_k(\omega)\right].
\end{equation*}
Consequently,
\begin{equation*}
	-\Delta_{\mathbb{S}^{d-1}} \psi_k=k(k+d-2) \psi_k.
\end{equation*}
The norm of $\psi_k$ can be computed explicitly. Consider
\begin{equation*}
	J_k =\int_{\mathbb{R}^d} e^{-|x|^2}(x_1^2+x_2^2)^k d x.
\end{equation*}
Write $x=r \omega$, with $r \in[0,\infty)$ and $\omega \in \mathbb{S}^{d-1}$. Since $|x_1+ix_2|^{2k}=(x_{1}^{2}+x_{2}^{2})^k$, $(x_{1}^2+x_{2}^{2})^k=r^{2k}(\omega_1^2+\omega_2^2)^k$ and $dx=r^{d-1}dr dS(\omega)$, one has
\begin{align*}
	J_k & =\int_0^{\infty} \int_{\mathbb{S}^{d-1}} r^{2 k}(\omega_1^2+\omega_2^2)^k e^{-r^2} r^{d-1} d S(\omega)d r   =\left(\int_0^{\infty} r^{2 k+d-1} e^{-r^2} d r\right)\left(\int_{\mathbb{S}^{d-1}}(\omega_1^2+\omega_2^2)^k d S(\omega)\right).
\end{align*}
Therefore,
\begin{equation*}
	J_k=\left(\int_0^{\infty} r^{2 k+d-1} e^{-r^2} d r\right)\Norm{\psi_k}{L^2\left(\mathbb{S}^{d-1}\right)}^2.
\end{equation*}
Using the change of variables $t=r^2$, one obtains
$$
J_k=\frac{1}{2} \Gamma(k+d/2)\Norm{\psi_k}{L^2(\mathbb{S}^{d-1})}^2.
$$
On the other hand, writing $\mathbb{R}^d=\mathbb{R}^2 \times \mathbb{R}^{d-2}$, using the change of variables $t=r^2$ once again and applying Fubini,
\begin{align*}
	J_k & =\left(\int_{\mathbb{R}^2} e^{-(x_{1}^2+x_{2}^{2})} (x_{1}^2+x_{2}^2)^{ k} d x_1dx_2\right)\left(\int_{\mathbb{R}^{d-2}} e^{-|x'|^2} dx'\right)  =\left(2 \pi \int_0^{\infty} e^{-r^2} r^{2 k+1} d r\right) \pi^{(d-2) / 2} \\
	& =\pi \Gamma(k+1) \pi^{(d-2) / 2}  =\pi^{d / 2} \Gamma(k+1)
\end{align*}

Comparing the two expressions for $J_k$ and using \cite[Eq. 5.11.12]{NewAbramowitz} gives
\begin{equation*}
	\Norm{\psi_k}{L^2\left(\mathbb{S}^{d-1}\right)}^2=2 \pi^{d / 2} \frac{\Gamma(k+1)}{\Gamma(k+d / 2)} \sim 2\pi^{d/2}k^{1-\frac{d}{2}}.
\end{equation*}

Suppose that $\supp \alpha \subset \omega=\{ x\in \mathbb{S}^{d-1}: x_{1}^2+x_{2}^{2}<e^{-4} \}$. The above asymptotics yield
\begin{equation*}
	\Norm{\alpha\phi_k}{L^2(\mathbb{S}^{d-1})}^2
	\lesssim \Norm{\phi_k}{L^2(\omega)}^{2}
	\leq C c_k^2 e^{-4k} \lesssim  k^{(d-2)/2}e^{-4k}
	\lesssim e^{-2k}.
\end{equation*}
Therefore, condition~\eqref{optimality: boundary general} holds with $C_0 = 2$.

\begin{figure}[htb!]
  \begin{center}
  \resizebox{7cm}{!}{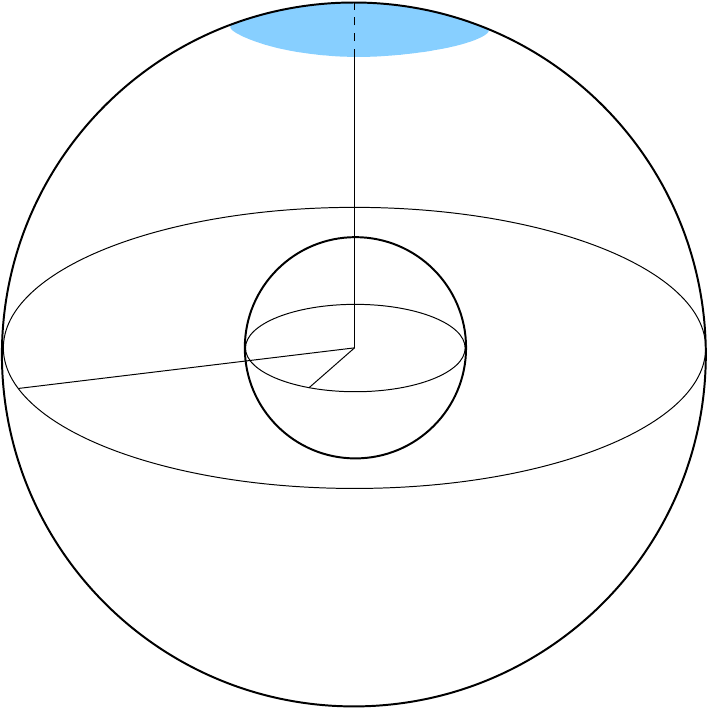}
  \caption{Sketch of an open set $\Omega\subset\mathbb{R}^3$ for which the hypotheses of Lemma~\ref{lem:quasimode-boundarygeneral} hold. (Some proportions are not respected.)}
  \label{fig:Ball}
\end{center}
\end{figure}

The example when $d=3$, corresponds to the manifold $\M = \mathbb S^2 \times [1/8,1]$ along with the metric $g = dy^2 + g_{\mathbb S^2}$, where $g_{\mathbb S^2}$ is the metric induced by the Euclidean metric on $\R^3$.
This corresponds to the case of the open set $\Omega$ of $\R^3$ given by 
$\Omega = B(0,1) \setminus \ovl{B(0,1/8)}$, where $B(0,r)$ denotes the Euclidean ball centered at $0$ with radius $r$, equipped with the non-Euclidean metric given by $g = dr^2 + g_{\mathbb S^2}$, in spherical coordinates. See Figure~\ref{fig:Ball} for a sketch.
\color{black}

\appendix 

\section{Well-posedness}

\begin{proof}[\bfseries Proof of Theorem \ref{wp-thm-1}]\label{proof:wp-thm-1}

Prove that the operator $(\SgGen, D(\SgGen))$, defined by~\eqref{wp-4}
and~\eqref{wp-6}, generates a $C_0$-semigroup of contractions
$e^{t\SgGen}$ on $\mathcal{H}$.  Suppose $U = (u, v, w) \in
D(\SgGen)$. A computation similar to those performed in
Proposition~\ref{proposition: spectrum} yields
\begin{equation*}
    \Re(\SgGen U, U)_{\mathcal{H}} 
    = -\|\alpha^{1/2} v|_{\Gamma_1}\|_{L^2(\Gamma_1)}^2 \leq 0,
\end{equation*}
that is, $\SgGen$ is dissipative. Moreover, $\SgGen^{-1} \colon
\mathcal{H} \to \mathcal{H}$ is bounded (even compact) with $D(\SgGen^{-1}) = R(\SgGen) = \mathcal{H}$; see the proof of
Proposition~\ref{proposition: spectrum}. Hence, $0$ lies in the resolvent set of $\SgGen$.
We claim that $R(\lambda \id_\H- \SgGen) = \H$ for  $\lambda > 0$ chosen \suff small.
Then, apply the Lumer-Phillips theorem and reach the conclusion; see \eg Theorem 4.3 in \cite[Chapter~1]{pazy}.
	
In fact, suppose $F\in \H$. For $\lambda>0$, define $G_\lambda:\H \to
\H$ given by $G U=\SgGen^{-1}(\lambda U-F)$. For $U, V \in \H$, one has
\begin{eqnarray*}
  \Norm{G_\lambda U-G_\lambda V}{\H}=\lambda \Norm{\SgGen^{-1}(U-V)}{\H}
  \leq \lambda \Norm{\SgGen^{-1}}{\mathcal{L}(\H,\H)} \Norm{U-V}{\H}.
\end{eqnarray*}
Hence, for $\lambda>0$ \suff small, $G_\lambda$ is a contraction map and has a
(unique) fixed point $U$ with $U \in D(A)$. Then, $\lambda U -
\SgGen U = F$.
\end{proof}

\section*{Acknowledgments}

R. Buffe was partially supported by the French National
Research Agency (ANR), Project QuBiCCS (ANR-24- CE40-3008). V. H. Gonzalez Martinez was supported by CAPES/COFECUB grant 88887.879175/2023-00 and CNPq grants 421573/2023-6, 402652/2024-0 and 300997/2025-6.

\bibliographystyle{acm}
\bibliography{references.bib}
\end{document}